\documentclass[10pt,twocolumn,twoside]{IEEEtran}

\usepackage{subfigure}
\usepackage{citesort}
\usepackage[fleqn]{amsmath}
\usepackage{ascmac}
\usepackage{citesort}
\usepackage[dvips]{graphicx}
\usepackage{psfrag}
\usepackage{amsmath}
\usepackage{amsbsy}
\usepackage{amssymb}
\usepackage{latexsym}
\usepackage{color}
\usepackage{theorem}

\graphicspath{{./images/}}

\newtheorem{proposition}{Proposition}
\newtheorem{example}{Example}
\newtheorem{algorithm}{Algorithm}
\newtheorem{definition}{Definition}
\newtheorem{theorem}{Theorem}
\newtheorem{lemma}{Lemma}

\newtheorem{fact}{Fact}
\newtheorem{remark}{Remark}

\newcommand{\migip}{\hfill $\blacksquare$}

\newcommand{\diag}{{\rm diag}}

\newcommand{\euclidspace}{{\mathcal{H}}}
\newcommand{\signal}[1]{{\boldsymbol{#1}}}
\newcommand{\Natural}{{\mathbb N}}
\newcommand{\norm}[1]{\left\|#1\right\|}
\newcommand{\abs}[1]{\left|#1\right|}

\newcommand{\real}{{\mathbb R}}

\newcommand{\innerprod}[2]{\left\langle{#1},{#2}\right\rangle}

\newcommand{\Fix}[1]{{\rm Fix}\left({#1}\right)}

\newcommand{\argmin}{\operatornamewithlimits{argmin}}

\newcommand{\gzh}{{\Gamma_0(\mathcal{H})}}
\newcommand{\prox}{{\rm s\mbox{-}Prox}}
\newcommand{\sprox}{{\rm s\mbox{-}Prox}}
\newcommand{\dom}{{\rm dom}~}
\newcommand{\range}{{\rm range}~}
\newcommand{\cont}{{\rm cont}~}
\newcommand{\interior}{{\rm int}~}

\definecolor{darkgreen}{rgb}{0,.6,0}
\definecolor{medorange}{rgb}{0.7,0.3,0}
\definecolor{cyancyan}{rgb}{0.68, 0.92, 0.92}

\usepackage{version}	
\begin{document}
%
%
%
\title{Monotone Lipschitz-Gradient Denoiser:
Explainability of Operator Regularization Approaches
Free From Lipschitz Constant Control}

%
%

 \author{Masahiro~{\sc Yukawa},~\IEEEmembership{Senior Member,~IEEE},
and Isao~{\sc Yamada},~\IEEEmembership{Fellow,~IEEE}
\thanks{This work was partially supported by JSPS Grants-in-Aid (23K22762).}
\thanks{M.~Yukawa is with the Department of Electronics and Electrical Engineering, Keio University,
Yokohama, Kanagawa 223-8522, Japan (e-mail: yukawa@elec.keio.ac.jp).}
\thanks{Isao Yamada is with 
the Department of Information and Communications Engineering, Institute of Science Tokyo, Meguro-ku, Tokyo 152-8550, Japan
(e-mail: isao@sp.ict.e.titech.ac.jp).}
\thanks{DOI: 10.1109/TSP.2025.3580667}
\thanks{© 2025 IEEE.  Personal use of this material is permitted.  Permission from IEEE must be obtained for all other uses, in any current or future media, including reprinting/republishing this material for advertising or promotional purposes, creating new collective works, for resale or redistribution to servers or lists, or reuse of any copyrighted component of this work in other works.}
}


%
%
%
 \markboth{IEEE TRANSACTIONS ON SIGNAL PROCESSING}
 {Yukawa and Yamada: Monotone Lipschitz-Gradient Denoiser}

%



\maketitle







\begin{abstract}

This paper addresses explainability of 
the operator-regularization approach
under the use of {\em monotone Lipschitz-gradient (MoL-Grad) denoiser} ---
an operator that can be expressed
as the Lipschitz continuous gradient of a differentiable convex function.
We prove that an operator is a MoL-Grad denoiser if and only if
it is the ``single-valued'' proximity operator of a weakly convex
function.
An extension of Moreau's decomposition is also shown
with respect to
a weakly convex function and the conjugate of its convexified function.
Under these arguments, two specific algorithms,
the forward-backward splitting algorithm and the primal-dual splitting
algorithm, are considered,
both employing MoL-Grad denoisers.
These algorithms generate a sequence of vectors converging
weakly, under conditions, to a minimizer of
a certain cost function which involves an ``implicit regularizer''
induced by the denoiser.
Unlike the previous studies of operator regularization, our framework
 requires no control of the Lipschitz constant in learning the denoiser.
The theoretical findings are supported by simulations.

\end{abstract}


%


\section{Introduction}\label{sec:intro}

Operator splitting (a.k.a.~proximal splitting) \cite{bauschke_book19}
is a standard technique nowadays for nonsmooth convex optimization.
It allows to take into account prior knowledge
such as sparsity and low-rankness
both of which are quite relevant to 
signal processing as well as machine learning.
Indeed, the naive discrete-valued measures (the $\ell_0$
pseudo-norm and the matrix rank) of 
such prior information are discontinuous, and
one often resorts to their convex relaxations
(the $\ell_1$ and nuclear norms)
which are typically nonsmooth.
Moreau's proximity operator resides
at the core of operator splitting algorithms,
and it has many examples ranging
from the soft-shrinkage and
the convex projection
to activation functions including
rectified linear unit (ReLU),
sigmoid, and softmax \cite{combettes_dnn20}, to name a few.
One of the remarkable advantages of the operator splitting algorithm is 
its computational efficiency as long as
the involved nonsmooth functions are prox friendly.
The celebrated examples of the operator splitting algorithm include
the proximal forward-backward splitting method (the proximal gradient
method) \cite{PFBS1,PFBS2},
the Douglas-Rachford splitting method \cite{PFBS1},
the alternating direction methods of multipliers \cite{glowinski75,gabay76},
and the primal-dual splitting methods \cite{chambolle11,condat13,vu13}.
See \cite{combettes18,condat_siam23}, 
for instance, for more about the history of the operator splitting algorithms.


Designing an objective function involving nonsmooth functions
is the first step 
in typical signal processing approaches based on operator splitting.
Then, an algorithm to minimize the objective
is either constructed or
selected from the available ones,
followed by a derivation of Moreau's proximity operators,
as well as other operators such as the gradient in most cases.
As the nonsmooth functions, the convex relaxations mentioned in the
previous paragraph are
usually chosen as a first choice.
Despite the mathematical tractability, such convex regularizers tend to
yield serious estimation biases \cite{antoniadis07} owing to the convexity.
One might therefore return to the first step to devise a better
objective function.
Indeed, nonconvex regularizers have been studied widely
\cite{fan01,chartrand08,candes08,zhang09,zhang,marjanovic12,shen18,yao18,Bwen18,chartrand07,soubies15,yukawa16,jeong14},
motivated by the fact that 
those regularizers are better approximations of the ideal $\ell_0$ pseudo-norm
and therefore reduce the estimation bias significantly.

Some promising approaches combining convex and nonconvex functions
have been studied \cite{dinh86,blake87,nikolova99}.
In these approaches, a weakly convex penalty is coupled with
a strongly convex fidelity so that 
the overall convexity of the whole objective function is preserved.
A generalized penalty function has been proposed in \cite{selesnick}
so that the overall convexity can be preserved even when applied to 
inverse problems of underdetermined linear systems
(where the fidelity function is not strongly convex).
In \cite{abe_ip20}, the penalty has further been generalized 
to permit a composition with a bounded linear operator,
and the generalized penalty has been shown to improve the performance of
the $\ell_1$-based total variation as a particular instance.
In \cite{yukawa23}, the convex-nonconvex strategy has been applied to 
outlier-robust signal recovery, and
a general framework has been established to clarify 
when the sufficient condition of overall convexity is also 
a necessary condition,
i.e., when the bound of the regularization parameter is
tight.
Its extension has also been presented in \cite{suzuki23}.


There are several other related approaches.
One is to directly devise a better ``operator''
that accommodates prior information 
in such an efficient way not to cause serious estimation
biases (see \cite{chartrand14} and the references therein).
Another popular one is an {\em operator-regularization} approach called
plug-and-play \cite{venkatakrishnan13},
where the proximity operator is replaced by another denoiser such as
deep neural networks.
``Explainability'' of the plug-and-play method
has been studied actively in recent years 
(see Section \ref{subsec:relation_to_prior_work}).
Most of the previous works in this line of research 
impose the strong assumption of (averaged) nonexpansiveness on the neural network, 
which often requires additional costs \cite{pesquet21,nair24}.
There are some works considering a weaker assumption 
allowing the Lipschitz constant to exceed one
\cite{hurault22,goujon24}.
Whereas all those works concern convergence of the iterates,
studies on optimality of the limit point are rather limited.
In the present work,
the following question is addressed
for operators which are not necessarily nonexpansive:
{\em when do the operator-regularization approaches actually minimize
some sort of objective function,
and how can such a function be characterized?}

In our approach, we suppose that an operator is given.
It could be designed by human, or it could also be learned by machine.
The designed/learned operator replaces the proximity operator
appearing in an operator splitting algorithm, and 
the optimization problem associated with the resultant algorithm
is ``derived'' subsequently.
As such, the processes of ``design'' and ``derivation''
are in the reverse order compared to the typical optimization-based approaches.
We clarify that the estimate
is not explicitly characterized 
as a minimizer (nor a stationary/critical point) of
a user-designed function.
It can only be characterized
as a minimizer of a convex function 
involving an {\em implicit regularizer} induced by a user-designed
``operator''.\footnote{
It can also be characterized as a fixed point of
some averaged nonexpansive operator from a wider viewpoint.}
A practical advantage of the proposed approach
is that the operator can be designed in a way that 
it can be computed simply but its corresponding cost function could
be complicated, or even impossible to express in a closed form.

To present a more concrete picture, let us consider the simple iterate:
\begin{equation}
 x_{k+1}:= T (x_k- \mu \nabla f(x_k)),~k\in\Natural,
\label{eq:pfbs}
\end{equation}
where $T$ is a nonlinear operator from a real Hilbert space $\euclidspace$ to
$\euclidspace$, 
$f:\euclidspace\rightarrow \real $ is a smooth convex function,
and $\mu\in\real_{++}$ is the step size.
Here, the term {\em smooth} (or {\em $\kappa$-smooth} more specifically) is
used in this paper to mean that the function is differentiable
with its gradient Lipschitz continuous with constant $\kappa>0$
(see \eqref{eq:def_lipschitz}) over $\euclidspace$.
If in particular $T$ is Moreau's proximity operator 
of a proper lower semicontinuous convex function $g$,
\eqref{eq:pfbs} is the classical forward-backward splitting algorithm
to minimize $\mu f + g$.
It will turn out that
$(x_k)_{k\in\Natural}$ converges
weakly to a minimizer (if exists) of a certain cost function
under the following condition (in addition to other technical assumptions).\footnote{
One may consider such a condition that
$T$ is a selection of the subdifferential of a (not necessarily
differentiable) proper lower-semicontinuous convex function $\psi$.
However, this seemingly weaker condition together with (Lipschitz)
continuity implies  Fr\'echet differentiability of $\psi$
\cite[Proposition 17.41]{combettes}.
To the best of authors' knowledge, the subdifferentiability is often
used in convex analysis to relate an operator to optimization.
The case of $\beta\geq 1$ corresponds to Moreau's proximity operator of
convex function \cite{moreau65},
 thus excluded in the present study.
}

\begin{definition}[MoL-Grad Denoiser]
\label{def:molgrad}
An operator $T:\euclidspace\rightarrow\euclidspace$
is said to be a monotone Lipschitz-gradient (MoL-Grad) denoiser
if it can be expressed as the gradient of 
a $\beta^{-1}$-smooth convex function for $\beta \in (0, 1)$;
i.e.,
$T= \nabla \psi$ for a Fr\'echet differentiable convex
function $\psi$ with
the $\beta^{-1}$-Lipschitz continuous gradient $\nabla \psi$.
\end{definition}
The term ``denoiser'' is used in Definition \ref{def:molgrad}
rather than the general term ``operator'',
because we primarily consider denoising operator (which is often
referred to as {\em denoiser} in the literature) for $T$.
The term ``monotone'' stems from the fact that
convexity of $\psi$ can be characterized by monotonicity of the gradient $\nabla \psi$
(see Section \ref{subsec:monotone}).
The convergence argument can be extended to 
other operator splitting algorithms than \eqref{eq:pfbs},
as elaborated later on.
It should be mentioned that
(i) the Lipschitz continuity imposed on MoL-Grad denoiser
is weaker than
(averaged) nonexpansiveness which requires the Lipschitz constant to be
one (or smaller),
(ii) every MoL-Grad denoiser is $\beta$-cocoercivity of $T$,
i.e., $\beta T$ is firmly nonexpansive, and
(iii) the present study includes nonseparable operators
in scope.\footnote{
For instance, when used 
as a penalty function for the linear inverse problem with
underdetermined systems,
a nonconvex function must be nonseparable to preserve
overall convexity;
i.e., there exists no separable nonconvex function that
makes the overall cost function convex \cite{selesnick_bayram_16}.
}
We mention that, when a given denoiser $T:\real^n\rightarrow \real^n$ is a
continuously differentiable mapping, it can be expressed as
$T=\nabla h$ for some $h:\real^n\rightarrow \real$
(i.e., $T$ is a conservative vector field)
if and only if $T$ possesses a symmetric Jacobian matrix at every point \cite{sreehari16}.

\subsection{Contributions}

In the first main part (Section \ref{sec:weak_convexity}), 
we link the operator $T$ to 
a generalized notion of Moreau's proximity operator for possibly nonconvex functions
to facilitate the analysis.
Here, the proximity operator of a nonconvex function is set-valued in general.
In the meanwhile, every MoL-Grad denoiser is continuous owing directly
to its Lipschitz continuity
(see Definition \ref{def:continuity} in Section \ref{subsec:nonexpansive});
if $T$ is discontinuous, the convergence analysis of the algorithm in
\eqref{eq:pfbs} would hardly be tractable in general.
As long as imposing continuity on $T$,
we can restrict ourselves to the {\em single-valued} case,
as seen from the following proposition.
\begin{proposition}
\label{proposition:single_valuedness}
Given a proper function $f:\euclidspace\rightarrow (-\infty,+\infty]$,
let $T:\euclidspace\rightarrow\euclidspace$ be a continuous\footnote{
Continuity of the denoiser is often assumed in the literature \cite{meinhardt17}.
} operator
such that
$T(x) \in \argmin_{y\in\euclidspace} [f(y) + (1/2)\norm{x-y}^2]\neq \varnothing$
for every $x\in\euclidspace$.
Then, the penalized function
$f + (1/2)\norm{x-\cdot}^2$
has a unique minimizer for every $x\in\euclidspace$.
\end{proposition}
\begin{proof}
See Appendix \ref{subsec:proof_prop1}.
\end{proof}

Proposition \ref{proposition:single_valuedness} states that,
if $f + (1/2)\norm{x-\cdot}^2$ has multiple minimizers,
there is no hope to find a continuous operator
$T:\euclidspace\rightarrow\euclidspace$
in the form of 
$T(x) \in \argmin_{y\in\euclidspace} [f(y) + (1/2)\norm{x-y}^2]$.
We therefore define our generalized proximity operator
as a unique minimizer of the penalized function only when the minimizer
exists uniquely.
Now, a natural question is what is a sufficient condition for 
$T$ to be continuous.
It is clear that 
$f + (1/2)\norm{x-\cdot}^2$ has a unique minimizer
(i.e., the necessary condition in the proposition holds true) 
if $f+ (1/2)\norm{\cdot}^2$ is strongly convex;
$f$ is weakly convex at least in this case.
This condition actually implies
that the denoiser $T$ can be expressed as the Lipschitz continuous gradient of a
convex function,
and thus $T$ is continuous.
Our findings of the first part are summarized below.
\begin{enumerate}
 \item An operator $T$ is a MoL-Grad denoiser
if and only if $T$ is the (generalized)
       proximity operator of a ($1-\beta$)-weakly convex function 
$\varphi = \psi^* -(1/2)\norm{\cdot}^2$
(Theorem \ref{theorem:weaklyconvex_necsuffcondition}).
This actually motivates the definition of 
MoL-Grad denoiser.

 \item 
 Moreau's decomposition is extended to weakly convex
       functions 
(Propositions \ref{proposition:extended_Moreau_decomposition}
and \ref{proposition:moreau_decomp_mod}).

\end{enumerate}

In the second part (Sections \ref{sec:applications}--\ref{sec:numerical}),
we consider the modified operator splitting algorithms with the
proximity operator replaced by a MoL-Grad denoiser,
analyzing the convergence based on the results established in the first part.
As a case study, we highlight two specific algorithms:
(i) the forward-backward splitting type algorithm given in
\eqref{eq:pfbs}, and 
(ii) the primal-dual splitting type algorithm which is basically 
the Condat--V\~u algorithm \cite{condat13,vu13}.
The specific contributions are listed below.
\begin{enumerate}
  \item[3)] 
The iterate \eqref{eq:pfbs} with a MoL-Grad denoiser $T$
converges weakly to a minimizer
	of $\mu f +  \varphi$ under conditions (Theorem
	  \ref{theorem:convergence}).
The monotonicity of our denoiser $T$ allows its associated regularizer
$\varphi$  to be (weakly) convex.

  \item[4)] 
 The primal-dual splitting type algorithm employing 
a MoL-Grad denoiser converges weakly to a minimizer
	of a certain function under conditions (Theorem \ref{theorem:primal_dual}).

\item[5)]  A systematic way of building operator-regularization
	algorithms using MoL-Grad denoisers with convergence guarantee
	is presented (Section \ref{subsec:systematic}).

\item[6)]
The proposed framework requires no control of the Lipschitz constant
during the learning process of the denoiser.
The implicit regularizer induced by the denoiser is not restricted to
	smooth ones (Remark \ref{remark:mol_grad_smoothness}).

\item[7)]  Simulation results support the theoretical findings
as well as showing potential advantages of the proposed approach
(Section \ref{sec:numerical}).
\end{enumerate}

A part of Sections \ref{sec:weak_convexity}, 
\ref{subsec:cocoercivity_optimization}, and \ref{subsec:examples}
has been presented in \cite{yukawa_sip23} without proofs.
The present work contains many new results/discussions
(such as Theorem \ref{theorem:primal_dual} and Propositions
\ref{proposition:single_valuedness} and
\ref{proposition:moreau_decomp_mod})
from a wider scope as well as proofs and simulations.

\section{Preliminaries}
\label{sec:preliminaries}

Let $(\euclidspace,\innerprod{\cdot}{\cdot})$ be a real Hilbert space
with the induced norm $\norm{\cdot}$.
Let $\real$, $\real_+$, $\real_{++}$, and $\Natural$ denote
the sets of real numbers, nonnegative real numbers,
strictly positive real numbers, and nonnegative integers, respectively.


\subsection{Properness, Subdifferentiability, and Lower Semicontinuity}
\label{subsec:properness}

We consider a function $f:\euclidspace\rightarrow (-\infty,+\infty]:=
\real\cup \{+\infty\}$.
A function $f$ is {\em proper} if the domain is nonempty; i.e.,
$\dom f:= \{x\in\euclidspace\mid f(x)<+\infty\}\neq \varnothing$.
Given a proper function $f$, 
the set 
\begin{equation}
~ \partial f(x):=\{z\in\euclidspace\mid \innerprod{y\! -\! x}{z} + f(x) \leq f(y), ~\forall y\in\euclidspace\}
\label{eq:def_subdifferential}
\end{equation}
is the subdifferential of $f$ at $x\in \euclidspace$ \cite{combettes};
each element $f'(x)\in \partial f(x)$ is a subgradient of $f$ at $x$.
A function $f$ is {\it lower semicontinuous} (or {\it closed}) on $\euclidspace$ 
if the level set
${\rm lev}_{\leq a} f:=
\left\{x\in\euclidspace: f(x)\leq a\right\}$
is closed for every $a\in\real$.
Every continuous function is lower semicontinuous.

\subsection{Convexity and Conjugation}
\label{subsec:convex}

A function $f:\euclidspace\rightarrow
(-\infty,+\infty]$ is convex on $\euclidspace$ if
$f(\alpha x + (1-\alpha)y)\leq 
\alpha f(x) + (1-\alpha)f(y)$ for every
$(x,y,\alpha)\in\dom f\times\dom f\times [0,1]$.
If the inequality of convex function holds with strict inequality
whenever $x\neq y$, $f$ is strictly convex.
For a positive constant $\rho\in\real_{++}$,
$f$ is $\rho$-strongly convex if $f- (\rho/2) \norm{\cdot}^2$ is convex,
while $f$ is $\rho$-weakly convex
if $f+ (\rho/2) \norm{\cdot}^2$ is convex.

Let $f:\euclidspace\rightarrow (-\infty,+\infty]$ be
a proper convex function.
Then, we have
(i) $\partial f(x)\neq \varnothing$ if $f$ is continuous at
$x\in\euclidspace$
\cite[Proposition 16.17]{combettes},
 and
(ii) $ \partial f(x)= \{\nabla f(x)\}$
if $f$ is G\^ateaux differentiable,
where $\nabla f(x)$ is the G\^ateaux gradient of $f$ at $x$.\footnote{
(G\^ateaux and Fr\'echet gradients of function)
Let  $U$ be an open subset of $\euclidspace$.
Then, a function $f:U\rightarrow \real$ is G\^ateaux differentiable 
at a point $x\in\euclidspace$ if there exists a vector $\nabla f(x)\in\euclidspace$
such that $\lim_{\delta\rightarrow 0} \frac{f(x+\delta h) -
f(x)}{\delta}= \innerprod{\nabla f(x)}{h}$ for every $h\in\euclidspace$.
On the other hand, 
a function $f:U\rightarrow \real$ is Fr\'echet differentiable 
at a point $x\in\euclidspace$
if there exists a vector $\nabla f(x)\in\euclidspace$
such that $\lim_{h\rightarrow 0} \frac{f(x+ h) -f(x) - \innerprod{\nabla
f(x)}{h}}{\norm{h}}= 0$.
In each case, $\nabla f(x)$ is the
G\^ateaux (or Fr\'echet) gradient  of $f$ at $x$.
If $f$ is Fr\'echet differentiable over $U$, it is also G\^ateaux
differentiable over $U$,
and the G\^ateaux and Fr\'echet gradients coincide.
}
Those readers who are unfamiliar with convex analysis may consider 
$\nabla f$ to be the standard gradient
without introducing themselves to the general notion of 
G\^ateaux differentiability.

The set of all proper lower-semicontinuous convex functions
$f:\euclidspace\rightarrow (-\infty,+\infty]$
is denoted by $\Gamma_0(\euclidspace)$.
Given a function $f\in\gzh$, 
{\em the Fenchel conjugate of $f$} is $f^*: \euclidspace \rightarrow
(-\infty,\infty]:x\mapsto
\sup_{y\in\euclidspace} \innerprod{x}{y} - f(y)$, satisfying (i)
$f^*\in\gzh$
and (ii) $u\in\partial f(x) \Leftrightarrow x \in\partial f^*(u)$.
Given a continuous convex function $f:\euclidspace\rightarrow \real$
and a positive constant $\rho\in\real_{++}$, the interplay between
smoothness and strong convexity is known \cite[Theorem 18.15]{combettes}:
$f$ is $\rho^{-1}$-smooth if and only if
the conjugate $f^*$ is $\rho$-strongly convex.

\subsection{Nonexpansiveness, Cocoercivity, and Fixed Point}
\label{subsec:nonexpansive}

In the present work, the following notions of convergence and continuity
will be used.

\begin{definition}[Strong/weak convergence]
\label{def:convergence}

Let $(x_k)_{k\in\Natural}\subset \euclidspace$  be a sequence
of vectors in a real Hilbert space $\euclidspace$.
 \begin{enumerate}
  \item The sequence $(x_k)$ converges strongly to a point
$\hat{x}\in\euclidspace$ if 
$\lim_{k\rightarrow \infty}\norm{x_k-x}=0$; 
in symbols, $x_k\rightarrow \hat{x}$ (as $k\rightarrow \infty$).

  \item The sequence $(x_k)$ converges weakly to a point
$\hat{x}\in\euclidspace$ if 
$\lim_{k\rightarrow \infty}\innerprod{x_k - \hat{x}}{y}=0$ for every
$y\in\euclidspace$; in symbols, $x_k\rightharpoonup \hat{x}$
	(as $k\rightarrow \infty$).
 \end{enumerate}
If $\euclidspace$ is finite dimensional, strong convergence coincides
 with weak convergence, i.e., weak convergence implies strong
 convergence,
and vice versa {\rm \cite[Fact 2.33]{combettes}}.
\end{definition}

\begin{definition}[Lipschitz continuity]
\label{def:continuity}

\begin{enumerate}
 \item 
An operator $T:\euclidspace\rightarrow \euclidspace$ is
Lipschitz continuous
with constant $\kappa>0$ (or $\kappa$-Lipschitz continuous for short)
if for every $x,y\in\euclidspace$
\begin{equation}
\norm{T(x)-T(y)}\leq \kappa \norm{x-y}.
\label{eq:def_lipschitz}
\end{equation} 
 \item 
An operator $T:\euclidspace\rightarrow \euclidspace$ is
continuous (as a mapping from the normed space
       $(\euclidspace,\norm{\cdot})$ to itself $(\euclidspace,\norm{\cdot})$)
if $T$ is continuous at every point $\hat{x}\in\euclidspace$,
i.e., for every sequence $(x_k)_{k\in\Natural}$,
$x_k\rightarrow \hat{x}$ implies $T(x_k)\rightarrow T(\hat{x})$.
\end{enumerate}
Lipschitz continuity of $T$ implies 
continuity of $T$, as can readily be verified 
by \eqref{eq:def_lipschitz}.
\end{definition}

A $1$-Lipschitz continuous operator is called {\em nonexpansive}.
Given a nonexpansive operator $N:\euclidspace\rightarrow \euclidspace$,
the operator $T:=(1-\alpha){\rm Id}+ \alpha N$ for $\alpha\in(0,1)$ is
{\em $\alpha$-averaged nonexpansive}.
In particular, (1/2)-averaged operator is called {\em firmly
nonexpansive}.
For $\beta\in\real_{++}$,
an operator $T:\euclidspace\rightarrow \euclidspace$ is {\em $\beta$-cocoercive}
if $\beta T$ is firmly nonexpansive.
A point $x\in\euclidspace$ is {\em a fixed point of $T$} if $T(x)=x$.
The set ${\rm Fix}(T):=\{x\in\euclidspace\mid T(x)=x\}$ of all fixed points
is called {\em the fixed point set of $T$}.

We define
the identity operator 
${\rm Id}:\euclidspace\rightarrow\euclidspace:x\mapsto x$
on $\euclidspace$.
The following fact is a part of
 \cite[Theorem 5.15]{combettes}
with simplification regarding the domain of the operator.

\begin{fact}
\label{fact:mann}
 Let $T:\euclidspace\rightarrow \euclidspace$ 
be a nonexpansive operator such that $\Fix{T}\neq \varnothing$.
Define $T_{\alpha_k}:= (1-\alpha_k) {\rm Id} + \alpha_k T$
for a real-number sequence $(\alpha_k)_{k\in\Natural}\subset [0,1]$ such that
$\sum_{k\in\Natural} \alpha_k(1-\alpha_k) = +\infty$.\footnote{
It suffices that $\alpha_k\in [\epsilon,1-\epsilon]$
($\forall k\in\Natural$) for some small $\epsilon\in(0,1/2)$.
}
Then, given an initial point $x_0\in\euclidspace$,
the sequence $(x_k)_{k\in\Natural}$ generated by
\begin{equation}
 x_{k+1} := T_{\alpha_k}(x_k),~~~k\in\Natural,
\label{eq:mann}
\end{equation}
converges weakly to a point $\hat{x}\in\Fix{T}$;
i.e.,
$x_k \rightharpoonup \hat{x}\in\Fix{T}$.

\end{fact}
The update rule in \eqref{eq:mann} is called {\it the Krasnosel'sk\u{i}-Mann
iterate}.
In particular, letting $\alpha_k:= \alpha$ $(\forall k\in\Natural)$ for
some $\alpha\in (0,1)$ in Fact \ref{fact:mann} implies that,
given any averaged nonexpansive mapping $T_\alpha$
with $\Fix{T_{\alpha}}\neq \varnothing$, 
it holds that
$T_\alpha^k (x_0) \rightharpoonup \hat{x}\in\Fix{T_{\alpha}}$.

\subsection{Monotonicity}
\label{subsec:monotone}

The subdifferential operator $\partial f:\euclidspace\rightarrow
2^\euclidspace$ defined in \eqref{eq:def_subdifferential}
maps a vector $x$ to the set (of subgradients).
In general, $T:\euclidspace\rightarrow 2^{\euclidspace}$ is called 
a set-valued operator,
where $2^{\euclidspace}$ is the power set of $\euclidspace$;
i.e., the family of all subsets of $\euclidspace$.
An operator $T:\euclidspace\rightarrow 2^{\euclidspace}$ is {\em monotone} if
\begin{equation}
\hspace*{-2em} \innerprod{x \! - \! y}{u \! - \! v}\geq 0,~\forall (x,u)\in{\rm gra~}T,~\forall
  (y,v)\in{\rm gra~}T,
\label{eq:def_monotone}
\end{equation}
where gra $T:=\{(x,u)\in\euclidspace^2\mid u\in T(x)\}$ is the graph of $T$.
A monotone operator $T:\euclidspace\rightarrow 2^{\euclidspace}$ is
{\em maximally monotone} if no other monotone operator
has its graph containing gra $T$ properly.
For instance, the subdifferential operator $\partial f$ of a proper function
$f:\euclidspace\rightarrow (-\infty,+\infty]$
is monotone, and 
it is maximally monotone if $f\in\Gamma_0(\euclidspace)$.
Suppose that a proper function
$f:\euclidspace\rightarrow (-\infty,+\infty]$ is 
G\^ateaux differentiable on $\dom f$ which is open and convex.
Then, $f$ is convex if and only if $\nabla f$ is monotone
\cite[Proposition 17.7]{combettes}.

\subsection{Proximity Operator and Moreau Envelope}
\label{subsec:prox}

Motivated by Proposition \ref{proposition:single_valuedness},
we define the single-valued proximity operator 
(the s-prox operator for short) of possibly nonconvex functions as below.

\begin{definition}[Single-valued proximity (s-prox) operator]
\label{def:sprox}

Let $f:\euclidspace\rightarrow (-\infty,+\infty]$
be a proper function.
Given a positive constant $\gamma > 0$,
assume that 
$f + (1/(2 \gamma)) \norm{x - \cdot}^2$ has a unique minimizer for every
 fixed $x\in\euclidspace$.
Then, the s-prox operator of $f$
of index $\gamma$ is defined by\footnote{
For proper (not necessarily convex) functions, the proximity operator is often defined 
as a set-valued operator
because the $f + (1/(2 \gamma)) \norm{x - \cdot}^2$ may have multiple
minimizers
\cite{gribonval20,bauschke21}
(see also https://proximity-operator.net/proximityoperator.html).
In the present work, we focus on the case in which
a unique minimizer exists over the entire space.
Although the multi-valued proximity operator coincides with our s-prox
operator when a unique minimizer exists, we give the specific name to
our operator and we use the different notation ${\sprox}_{\gamma f}$ to
convey our message that there are remarkable advantages from 
the viewpoints of explainability and optimization if we restrict ourselves
to the unique-minimizer case.
}
\begin{equation}
\hspace*{-1.3em}
{\sprox}_{\gamma f}:
\euclidspace \! \rightarrow \! \euclidspace:
x \mapsto \argmin_{y \in \mathcal{H}}
\!\Big(f (y) + \frac{1}{2 \gamma}
\norm{x \!- \!y}^2\Big).
\label{eq:def_prox}
\end{equation}
 
\end{definition}

In the particular case of $f\in\gzh$,
existence and uniqueness of minimizer is automatically
ensured owing to the strict convexity of $\norm{\cdot}^2$
and the coercivity of $\norm{\cdot}^2$
(i.e., $\norm{y}^2\rightarrow +\infty$ as
$\norm{y}\rightarrow +\infty$).
In this convex case, 
${\sprox}_{\gamma f}$ reduces to the classical Moreau's proximity operator
\cite{moreau62,moreau65,combettes05}
 which is firmly nonexpansive,
and the following identity (known as {\em Moreau's decomposition}) holds:
\begin{equation}
 {\rm Id}= \prox_{\gamma f} + \gamma
 \prox_{\gamma^{-1} f^*}\circ (\gamma^{-1}{\rm Id}).
\label{eq:classical_moreau_decomposition}
\end{equation}

Given a function $f\in\Gamma_0(\euclidspace)$,
its Moreau envelope \cite{moreau62,moreau65,combettes05} of index $\gamma\in(0,\infty)$ is defined
as follows:
\begin{equation}
 ^{\gamma} f:\euclidspace\rightarrow\real:
 x\mapsto \min_{y\in\euclidspace}
\Big(
f(y) + \frac{1}{2\gamma}\norm{x-y}^2
\Big).
\end{equation}

The Moreau envelope $^{\gamma}f$
is convex and Fr\'echet differentiable
with gradient given by
\begin{equation}
\nabla ~ ^{\gamma}f (x)= \gamma^{-1}
\left(x - \prox_{\gamma f} (x)\right),
\end{equation}
which is $\gamma^{-1}$-Lipschitz continuous
\cite{moreau62,moreau65,combettes05,yyy_springer_book11}.
This means that
the Moreau envelope is a smooth approximation of
a potentially nondifferentiable convex function $f\in\gzh$.
See \cite{yyy_springer_book11} for more details.



\section{Fundamental Results}\label{sec:weak_convexity}

The operator $\sprox_{\gamma f}$ in \eqref{eq:def_prox}
is well-defined if $\gamma f$ is
$\rho$-weakly convex for $\rho\in(0,1)$
so that $\gamma f + (1/2)\norm{x - \cdot}^2$ is strongly (and thus strictly)
convex and coercive.
This simple observation is thoroughly investigated,
followed by an extension of Moreau's decomposition.

\subsection{Weakly Convex Function and
Gradient Operator of Smooth Convex Function}
\label{subsec:wc_and_coco}

We show that a given nonlinear operator
$T:\euclidspace\rightarrow\euclidspace$
is the proximity operator of some weakly convex function
if and only if it can be expressed as the gradient operator of a smooth convex function,
i.e., a MoL-Grad denoiser (see Definition \ref{def:molgrad}).
The proof relies on the following lemma.

\begin{lemma}
\label{lemma:prox_resolvent_nonconvex_case}
For a proper function $h:\euclidspace\rightarrow (-\infty,+\infty]$
and a vector $x\in\euclidspace$, 
let $J_x(y):= h(y)+ (1/2)\norm{y-x}^2$, $y\in\euclidspace$.
Then, the following equalities hold:
\begin{align}
 \partial J_x(y) = &~ \partial \Big(h + \frac{1}{2}\norm{\cdot}^2 \Big) (y) -x, \label{eq:subdifferential_decomposition}\\
\argmin_{y\in\euclidspace} J_x(y)  =&~ \Big[\partial \Big(h +
 \frac{1}{2}\norm{\cdot}^2\Big) \Big]^{-1}(x)
\label{eq:prox_resolvent_nonconvex}\\
:=&~ \Big\{p\in\euclidspace\mid x \in\partial \Big(h +
 \frac{1}{2}\norm{\cdot}^2\Big)(p) \Big\}.
\nonumber
\end{align}
In particular, if $h + (1/2)\norm{\cdot}^2$ is coercive and strictly convex,
$J_x$ has a unique minimizer, which is
\begin{equation}
 \sprox_h(x) = \Big[\partial \Big(h +  \frac{1}{2}\norm{\cdot}^2\Big)\Big]^{-1}(x)
\label{eq:prox_subdiff}
\end{equation}
(in the sense of \eqref{eq:def_prox}).

\end{lemma}

\begin{proof}
Fix $y\in\euclidspace$ arbitrarily.
For $s\in\euclidspace$, the following equivalence holds:
\begin{align}
\hspace*{-2em}& ~s \in \partial J_x(y)
\nonumber\\
\hspace*{-2em}\Leftrightarrow  & ~
\innerprod{z-y} {s} + J_x(y)  \leq 
J_x(z), ~\forall z\in\euclidspace
 \nonumber\\
\hspace*{-2em}\Leftrightarrow & ~ \innerprod{z-y} {s} + h(y) +  \frac{1}{2}\norm{y}^2 -\innerprod{y}{x}
\nonumber\\
\hspace*{-2em} &
\hspace*{7.5em} \leq h(z) + \frac{1}{2}\norm{z}^2 -\innerprod{z}{x}, ~\forall z\in\euclidspace
\nonumber\\
\hspace*{-2em}\Leftrightarrow & 
 \innerprod{z  -y} {s+x} + h(y) + \frac{1}{2}  \norm{y}^2 
 \leq   h(z)
+  \frac{1}{2}  \norm{z}^2 \!\!\!, ~ \forall z \in \euclidspace
\nonumber\\
%
%
\hspace*{-2em}\Leftrightarrow & ~
  s\in \partial \Big(h + \frac{1}{2}\norm{\cdot}^2\Big)(y)-x,
\nonumber
\end{align}
which verifies \eqref{eq:subdifferential_decomposition}.
Now, let $p\in \euclidspace$.
Then, from \eqref{eq:subdifferential_decomposition}, the following
 equivalence holds:
\begin{align*}
\hspace*{-1.6em}
 p\in \argmin_{y\in\euclidspace} J_x(y)
\Leftrightarrow 
0\in \partial J_x(p) = \partial \Big(h +
 \frac{1}{2}\norm{\cdot}^2\Big)(p)-x &\\
\hspace*{-1.7em}
\Leftrightarrow
x\in \partial \Big(h + \frac{1}{2}\norm{\cdot}^2\Big)(p)
\Leftrightarrow 
p\in \Big[\partial \Big(h + \frac{1}{2}\norm{\cdot}^2\Big)\Big]^{-1}(x)&,
\end{align*}
which verifies \eqref{eq:prox_resolvent_nonconvex}.

When $h + (1/2)\norm{\cdot}^2$ is coercive and strictly convex,
the coercivity and strict convexity of $J_x$ ensures
the existence and uniqueness of its minimizer, respectively,
and thus \eqref{eq:prox_resolvent_nonconvex}, together with the definition of
the proximity operator in \eqref{eq:def_prox}, implies \eqref{eq:prox_subdiff}.
\end{proof}

\begin{theorem}
\label{theorem:weaklyconvex_necsuffcondition}

 Let $T:\euclidspace\rightarrow \euclidspace$.
Then, for every $\beta\in(0,1)$, the following two conditions are equivalent.
\begin{enumerate}
 \item[(C1)] $T=\sprox_\varphi$
       for some $\varphi:\euclidspace\rightarrow (-\infty,+\infty]$
such that $\varphi + \left((1-\beta)/2\right)\norm{\cdot}^2\in\Gamma_0(\euclidspace)$.
 \item[(C2)] $T$ is a MoL-Grad denoiser (see Definition \ref{def:molgrad}), i.e.,
the following hold jointly.
\begin{enumerate}
 \item [1)] $T= \nabla \psi$ for some Fr\'echet differentiable convex function
       $\psi\in \gzh$.

 \item [2)] $T$ is $\beta$-cocoercive, or equivalently
$\beta^{-1}$-Lipschitz continuous (by the Baillon Haddad theorem {\rm \cite{baillon77}}).
\end{enumerate}
\end{enumerate}

In particular, the following statements hold.
\begin{enumerate}
 \item[(a)]
Assume that {\rm (C1)} is satisfied.
Define $\check{\varphi}:=\varphi + (\left(1-\beta\right)/2) \norm{\cdot}^2\in\Gamma_0(\euclidspace)$.
Then, it holds that
\begin{equation}
T= \sprox_\varphi = \nabla \underbrace{\Big(\varphi + \frac{1}{2}\norm{\cdot}^2\Big)^*}_{=\psi}=
\nabla ~^\beta(\check{\varphi}^*),
\end{equation}
which is $\beta$-cocoercive (thus $\beta^{-1}$-Lipschitz continuous and
	   maximally monotone).

 \item[(b)] 
Assume that {\rm (C2)} is satisfied.
Then, it holds that
\begin{equation}
T= \nabla \psi = \sprox_{\psi^* -(1/2)\norm{\cdot}^2},
\end{equation}
where 
\begin{equation}
\varphi = \psi^* - \frac{1}{2} \norm{\cdot}^2
\end{equation}
is $(1-\beta)$-weakly convex.
\end{enumerate}

\end{theorem}
\begin{proof}
 C1 $\Rightarrow$ C2:
For the sake of readability, we present a proof which the reader may
follow directly with a minimal use of the known results in convex
 analysis.\footnote{
The proof can be shortened by using \cite[Example 22.7]{combettes} and 
\cite[Example 14.1]{combettes} to verify 
\eqref{eq:cocoercivity_of_T} and
\eqref{eq:conjugate_of_sum}.
}
Owing to the convexity of $\check{\varphi}:=\varphi + ((1-\beta)/2)\norm{\cdot}^2\in \Gamma_0(\euclidspace)$,
we have
$T^{-1} = \sprox_{\varphi}^{-1} 
= \partial (\varphi+(1/2)\norm{\cdot}^2) =  
\partial (\check{\varphi} + (\beta/2)\norm{\cdot}^2) = 
\partial \check{\varphi}  + \beta {\rm Id}$
by Lemma  \ref{lemma:prox_resolvent_nonconvex_case}.
Given any fixed $x,y\in\euclidspace$,
let $u\in T^{-1}(x)$ and $v\in T^{-1}(y)$
($\Leftrightarrow x=T(u),~y=T(v)$).
Then, since $T^{-1} - \beta {\rm Id}= \partial \check{\varphi}$ is a monotone operator,
it follows that
\begin{align}
\hspace*{-1em}& \innerprod{x-y}{(u-\beta x) - (v-\beta y)}\geq 0 \nonumber\\
\hspace*{-1em}\Leftrightarrow & \innerprod{x-y}{u - v}\geq \beta \norm{x-y}^2 \nonumber\\
\hspace*{-1em}\Leftrightarrow & \innerprod{T(u)-T(v)}{u - v}\geq \beta \norm{T(u) -
 T(v)}^2 \nonumber\\
\hspace*{-1em}\Leftrightarrow & \innerprod{\beta T(u)-\beta T(v)}{u - v}\geq \norm{\beta T(u) -
 \beta T(v)}^2,
\label{eq:cocoercivity_of_T}
\end{align}
where the last inequality implies that $T$ is $\beta$-cocoercive.
On the other hand, 
for every $y\in\euclidspace$,
it holds that
\begin{align}
\hspace*{-2em}
[^\beta(\check{\varphi}^*)]^*(y) =& \sup_{x\in\euclidspace}
\Big[
\innerprod{x}{y} - \inf_{z\in\euclidspace} \Big(
\check{\varphi}^*(z) + \frac{1}{2\beta}\norm{z-x}^2
\Big)
\Big]
\nonumber\\
=& \sup_{x\in\euclidspace}
\Big[\!
\innerprod{x}{y} + \sup_{z\in\euclidspace} \Big(
\!\! -\check{\varphi}^*(z) - \frac{1}{2\beta}\norm{z-x}^2
\Big)
\Big]
\nonumber\\
=& \sup_{z\in\euclidspace}
\Big[\!
-\check{\varphi}^*(z)
+ \sup_{x\in\euclidspace} \Big(
\! \innerprod{x}{y} 
 - \frac{1}{2\beta}\norm{z-x}^2
\!\Big)
\Big]
\nonumber\\
=& \sup_{z\in\euclidspace}
\Big[
-\check{\varphi}^*(z)
+ 
\innerprod{z}{y} +\frac{\beta}{2}\norm{y}^2
\Big]
\nonumber\\
=&~ 
\check{\varphi}^{**}(y)
 +\frac{\beta}{2}\norm{y}^2
=
\check{\varphi}(y)
 +\frac{\beta}{2}\norm{y}^2,
\label{eq:moreau_envelop_conjugate}
\end{align}
where $\check{\varphi}^{**}=\check{\varphi}$ because $\check{\varphi}\in\Gamma_0(\euclidspace)$
by virtue of the Fenchel-Moreau theorem \cite[Theorem 13.37]{combettes}.
From \eqref{eq:moreau_envelop_conjugate}, it follows that
\begin{equation}
\Big(
\check{\varphi} +\frac{\beta}{2}\norm{\cdot}^2
\Big)^*
=
[^\beta(\check{\varphi}^*)]^{**}
=
~ ^\beta(\check{\varphi}^*) 
=:
\psi,
\label{eq:conjugate_of_sum}
\end{equation}
where 
$^\beta(\check{\varphi}^*)$,
and thus $\psi$, is
a Fr\'echet differentiable function with
$\beta^{-1}$-Lipschitz continuous gradient 
(see Section \ref{subsec:prox}):
\begin{equation}
\nabla \psi
=
\nabla
\Big(
\check{\varphi} +\frac{\beta}{2}\norm{\cdot}^2
\Big)^*
=\Big[
\partial \Big(
\check{\varphi} +\frac{\beta}{2}\norm{\cdot}^2
\Big)
\Big]^{-1} = T.
\end{equation}
The proof in this part also verifies Theorem \ref{theorem:weaklyconvex_necsuffcondition}(a).

\noindent C2 $\Rightarrow$ C1:
We observe that
$T^{-1}= (\nabla \psi)^{-1}=\partial \psi^*
 = \partial (\varphi+(1/2)\norm{\cdot}^2)$, 
where $\varphi:= \psi^*-(1/2)\norm{\cdot}^2$.
This leads to 
$T=[\partial (\varphi + (1/2)\norm{\cdot}^2)]^{-1} = \sprox_{\varphi}$.
Here, the Fr\'echet differentiability of $\psi$ and the
$\beta^{-1}$-Lipschitz continuity of $\nabla \psi$
imply $\beta$-strong convexity of the conjugate $\psi^*$
(see Section \ref{subsec:convex}), and
we thus obtain 
$\varphi + \left((1-\beta)/2\right)\norm{\cdot}^2
= \psi^*- (\beta/2)\norm{\cdot}^2  \in\Gamma_0(\euclidspace)$.
The proof in this part also verifies Theorem \ref{theorem:weaklyconvex_necsuffcondition}(b).
\end{proof}

\begin{figure}[t!]
\centering
 \includegraphics[width=8.5cm]{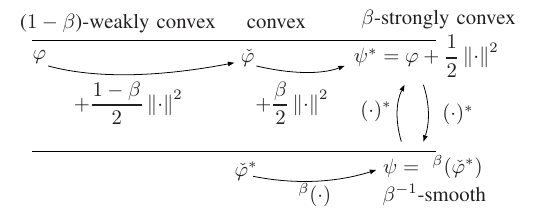}
 \caption{A weakly convex function $\varphi$, 
a smooth function $\psi$, and the intermediate ``neutrally'' convex function $\check{\varphi}$ between
$\varphi$ and the strongly convex  function $\psi^*$.}
\vspace{0em}
\label{fig:functions}
\end{figure}

The relations among $\varphi$, $\psi$, and $\check{\varphi}$ are depicted in
Fig.~\ref{fig:functions}.
An application of 
Theorem \ref{theorem:weaklyconvex_necsuffcondition}(b)
has been shown in \cite{suzuki_icassp24}, where
an external division of two Moreau's proximity operators is shown to be
the proximity operator of a weakly convex function under a certain
condition (see Example \ref{example:doscar} in 
Section \ref{subsec:shrinkage_structured}).

\begin{remark}[Alternative proof of Theorem
 \ref{theorem:weaklyconvex_necsuffcondition}(b)]

Theorem \ref{theorem:weaklyconvex_necsuffcondition}(b)
could be proved using known results as follows.
Since $\nabla \psi$ is $\beta^{-1}$-Lipschitz continuous
under condition (C2), 
it follows by {\rm \cite[Theorem 18.15]{combettes}}
that $\psi={}^\beta(\check{\varphi}^*)$,
where $\check{\varphi}:= \psi^* -
 (\beta/2)\norm{\cdot}^2\in\Gamma_0(\euclidspace)$.
Hence, by Moreau's decomposition, it follows that
$ \nabla \psi = \beta^{-1}({\rm Id}- \prox_{\beta \check{\varphi}^*})
= \prox_{\beta^{-1}\check{\varphi}}\circ (\beta^{-1}{\rm Id})
= \prox_{\beta^{-1}\psi^* - (1/2)\norm{\cdot}^2}\circ (\beta^{-1}{\rm Id}) 
= \sprox_{\psi^* -(1/2)\norm{\cdot}^2},
$
where the last inequality can be verified by
the definition of proximity operator with some manipulations.
Here, it is clear that
$\varphi:=\psi^* -(1/2)\norm{\cdot}^2
= \check{\varphi} - ((1-\beta)/2) \norm{\cdot}^2$ is
$(1-\beta)$-weakly convex.
\end{remark}

\subsection{Extension of Moreau's Decomposition}

We show that Moreau's decomposition (see Section \ref{subsec:prox})
can be extended to weakly convex functions.

\begin{proposition}
\label{proposition:extended_Moreau_decomposition}
Let 
$\varphi:\euclidspace\rightarrow (-\infty,+\infty]$ be a proper function
 such that 
$\check{\varphi}:=\varphi + ((1-\beta)/2) \norm{\cdot}^2\in\Gamma_0(\euclidspace)$
for $\beta\in(0,1)$.
Then, the following equalities hold:
\begin{align}
  \sprox_{\varphi} & = \prox_{\check{\varphi}/\beta} \circ
  \beta^{-1} {\rm Id},
\label{eq:proxphi_proxf}\\
{\rm Id}& = \beta \sprox_{\varphi} + \prox_{\beta \check{\varphi}^*}.
\label{eq:generalized_moreau_identity}
\end{align}
\end{proposition}

\begin{proof}
By using Moreau's decomposition,
it can be verified that
\begin{align}
& 
\beta^{-1}\left( {\rm Id}- \prox_{\beta \check{\varphi}^*} \right)(x)
 = \prox_{\beta^{-1} \check{\varphi}} (\beta^{-1}x)
\nonumber\\
& = \argmin_{y\in\euclidspace} 
\left(
 \beta^{-1}  \check{\varphi}(y) + \frac{1}{2}
\norm{y - \beta^{-1}x }^2
\right)
\nonumber\\
& = \argmin_{y\in\euclidspace} 
\left(
\varphi (y) +
\frac{1-\beta}{2} \norm{y}^2 + \frac{\beta}{2} \norm{y-\beta^{-1}x }^2
\right)
\nonumber\\
& = \argmin_{y\in\euclidspace} 
\left(
\varphi (y) +
\frac{1}{2} \norm{y-x}^2
+ K_x
\right)
= \sprox_{\varphi} (x),
\nonumber
\end{align}
where $K_x := ((\beta^{-1} -1)/2) \norm{x}^2$ is constant in $y$.
\end{proof}

The identity in \eqref{eq:proxphi_proxf} can be written equivalently as follows:
\begin{equation}
\underbrace{\prox_{\check{\varphi}/\beta}}_{\mbox{firmly nonexpansive}}
= \underbrace{\sprox_{\varphi}}_{\beta\mbox{-cocoercive}}
\circ \beta {\rm Id}.
\label{eq:moreau_decomposition_equivalent}
\end{equation}
The relation in \eqref{eq:moreau_decomposition_equivalent} 
implies $\beta$-cocoercivity of $\sprox_{\varphi}$ 
(see Lemma \ref{lemma:cocoercivity} in
Appendix \ref{appendix:averagedness_preservation}), which
is consistent
with Theorem \ref{theorem:weaklyconvex_necsuffcondition}.

\begin{remark}
[On Proposition \ref{proposition:extended_Moreau_decomposition}]
As opposed to the case of $\sprox_{\varphi}$ for the weakly convex
 function  $\varphi$,
$\prox_{\beta^{-1} \check{\varphi}}$ and $\prox_{\beta \check{\varphi}^*}$ are 
the proximity operators in the sense of Moreau,
because $\beta^{-1} \check{\varphi}$ and $\beta \check{\varphi}^*$ are convex.
The identity in \eqref{eq:proxphi_proxf}  gives an expression of 
the s-prox operator of the weakly convex function in terms of
the ``classical'' proximity operator of the convex function $\beta^{-1}\check{\varphi}$.
While the identity in \eqref{eq:generalized_moreau_identity} resembles 
Moreau's decomposition, the essential difference is that $\varphi$ is
 weakly convex. 
More specifically, \eqref{eq:generalized_moreau_identity} gives a way of
 decomposing a vector using the proximity operators of the weakly convex
function $\varphi$ and the convex function $\check{\varphi}^*$.
\end{remark}

Proposition \ref{proposition:extended_Moreau_decomposition} 
can be given in a general form as below.

\begin{proposition}
\label{proposition:moreau_decomp_mod}

Let $\phi:\euclidspace\rightarrow (-\infty,+\infty]$
such that $\check{\phi}:=\phi +
 (\rho/2)\norm{\cdot}^2\in\Gamma_0(\euclidspace)$
for a positive constant $\rho\in\real_{++}$.
Then, for every positive constant $\sigma\in\real_{++}$, it holds that

\begin{align}
 \hspace*{-1.5em}
& \prox_{\check{\phi}/\sigma} =
\sprox_{\phi/(\sigma+\rho)} \circ \left(
\frac{\sigma}{\sigma+\rho}
{\rm Id}
\right),
\label{eq:prox_relations_in_lemma_proof}\\
&  {\rm Id} = 
 \prox_{\sigma \check{\phi}^*} +
 \sigma \sprox_{\phi/(\sigma+\rho)} \circ (\sigma+ \rho)^{-1}{\rm Id}.
\end{align}
\end{proposition}

\begin{proof}
Let $\varphi:= \phi/(\sigma+\rho)$ in \eqref{eq:moreau_decomposition_equivalent}.
Then, $\varphi$ is $(1-\beta)$-weakly convex for $\beta:= \sigma/(\sigma+\rho)\in(0,1)$,
and we thus obtain \eqref{eq:prox_relations_in_lemma_proof}.
Using Moreau's decomposition
\eqref{eq:classical_moreau_decomposition} and then applying
\eqref{eq:prox_relations_in_lemma_proof} yield
\begin{align}
 {\rm Id}  =&~ \prox_{\sigma \check{\phi}^*} + \sigma \prox_{\check{\phi}/\sigma}\circ
 (\sigma^{-1}{\rm Id}) \nonumber\\
=&~ 
\prox_{\sigma \check{\phi}^*} +
\sigma \sprox_{\phi/(\sigma+\rho)} 
\circ 
   (\sigma+\rho)^{-1}
   {\rm Id}.
\end{align}
\end{proof}

 



\section{Operator-Regularization Approaches --- 
Characterization of Implicit Regularizer and
Convergence to Optimal Point}
\label{sec:applications}

We have seen that a MoL-Grad denoiser
can always be expressed as
$T=\sprox_\varphi$ for the $(1-\beta)$-weakly convex function
$\varphi= \psi^* - (1/2)\norm{\cdot}^2$
(Theorem \ref{theorem:weaklyconvex_necsuffcondition}).
We use this fundamental result to make the plug-and-play methods 
``transparent (a while box)''.
More specifically, we consider operator splitting algorithms
with the proximity operators replaced by
our MoL-Grad denoiser $T$,
which accommodates prior information (often implicitly).

We shall show that such plug-and-play methods 
(employing our denoiser $T=\nabla \psi$)
asymptotically optimize a certain cost function
in the sense of generating a vector sequence 
convergent to its minimizer.
Here, the cost function involves an {\em implicit ``operator-induced'' regularizer}
which depends on
$\varphi = \psi^* - (1/2)\norm{\cdot}^2$.
In this way, the loss function is implicitly regularized
by the denoiser $T$, as opposed to the case of
the traditional {\em functional-regularization} approaches 
which explicitly regularize the loss by some penalty functions.

In the following, two specific {\em operator-regularization} algorithms
are studied.
A generic way of applying the same idea systematically to different algorithms
is then presented.

\subsection{Forward-Backward Splitting Type Algorithm}
\label{subsec:cocoercivity_optimization}

Based on the results given in Section \ref{sec:weak_convexity},
we show that the algorithm in \eqref{eq:pfbs} with a MoL-Grad denoiser $T$
converges to a solution of an optimization problem.

\begin{theorem}
\label{theorem:convergence}
Let $f\in \gzh$ be a $\kappa$-smooth $\rho$-strongly-convex function with
$\kappa>\rho>0$.
Assume that $T:=\nabla \psi:\euclidspace \rightarrow \euclidspace$
is a MoL-Grad denoiser (see Definition \ref{def:molgrad})
for
 $\beta\in((\kappa-\rho)/(\kappa+\rho),1)\subsetneq(0, 1)$
so that 
the Lipschitz constant is bounded by
$\beta^{-1}< (\kappa+\rho)/(\kappa-\rho)$.
Set
$\mu\in [(1-\beta)/\rho, (1+\beta)/\kappa)$.
Let $\varphi:=\psi^* - (1/2)\norm{\cdot}^2$.
Then, the following hold.
\begin{enumerate}
 \item Let
$\hat{f} := f- [(1-\beta)/(2\mu)] \norm{\cdot}^2 \in \Gamma_0(\euclidspace)$ and
$\check{\varphi} := \varphi + [(1-\beta)/2] \norm{\cdot}^2 \in \Gamma_0(\euclidspace)$
so that $\mu \hat{f} + \check{\varphi} =\mu f + \varphi$.
Then, it holds that
\begin{equation}
T\circ ({\rm Id}- \mu \nabla f) = 
\prox_{\beta^{-1}\check{\varphi}}\circ ({\rm Id}- \beta^{-1}\mu \nabla \hat{f})
\end{equation}
with $ \beta^{-1} \in (1,2/L_{\mu \nabla \hat{f}})$, where
$L_{\mu \nabla \hat{f}}:= \mu\kappa - (1-\beta)>0$ is the Lipschitz
       constant of $\mu \nabla \hat{f}$.

 \item Suppose that 
$\mu f + \varphi$
has a minimizer in $\euclidspace$.
Then, for an arbitrary $x_0\in\euclidspace$,
 the sequence $(x_k)_{k\in\Natural} \subset\euclidspace$ generated by  \eqref{eq:pfbs} 
converges weakly to a minimizer $\hat{x}\in\euclidspace$ of $\mu f + \varphi$.
(In this case, $\varphi$ is the implicit regularizer.)
\end{enumerate}

\end{theorem}

\begin{proof}
1) Since $\mu\geq (1-\beta)/\rho \Leftrightarrow \rho \geq
 (1-\beta)/\mu$, 
the $\rho$-strong convexity of $f$ immediately implies
 $\hat{f}\in\Gamma_0(\euclidspace)$.
On the other hand, since
$\varphi$ is $(1-\beta)$-weakly convex by 
Theorem \ref{theorem:weaklyconvex_necsuffcondition},
$\hat{\varphi}\in\Gamma_0(\euclidspace)$ follows.
Now, using Theorem \ref{theorem:weaklyconvex_necsuffcondition}(b) and
Proposition \ref{proposition:extended_Moreau_decomposition}, we can
verify that
\begin{align*}
\hspace*{-1em} T \circ ({\rm Id}- \mu \nabla f) & = \sprox_{\varphi}\circ ({\rm Id}- \mu \nabla f)
\\
& = \prox_{\beta^{-1}\check{\varphi}}\circ (\beta^{-1}{\rm Id}- \beta^{-1}\mu \nabla f)\\
& = 
\prox_{\beta^{-1}\check{\varphi}}\circ ({\rm Id}- \beta^{-1}\mu 
\nabla \hat{f}).
\end{align*}
It can be verified that
$L_{\mu \nabla \hat{f}}
= \mu\kappa-(1-\beta)
> \mu\rho-(1-\beta)\geq 0$, 
where the assumption
$\mu \geq (1-\beta)/\rho$ is used again
to verify the last inequality.
On the other hand, we observe that
$\mu<(1+\beta)/\kappa
  \Leftrightarrow
\mu\kappa- (1 - \beta) < 2\beta
  \Leftrightarrow 
\beta^{-1} < 2/(\mu\kappa - (1-\beta)) = 2/L_{\mu \nabla \hat{f}}$,
and thus $\beta^{-1}<2/L_{\mu \nabla \hat{f}}$ is verified.
Finally, since $\nabla f$ is $\kappa$-Lipschitz continuous by the assumption,
it is not difficult to see that
$\mu\nabla \hat{f} = \mu \nabla f - (1-\beta){\rm Id}$ is 
$(\mu\kappa- (1-\beta))$-Lipschitz continuous
(see, e.g., \cite[Lemma 8]{bayram16}).

\noindent 2)
By Theorem \ref{theorem:convergence}.1, 
\eqref{eq:pfbs} can be regarded as
the forward-backward splitting iterate for the sum of the convex
functions $\mu \hat{f}$ and $\check{\varphi}$ with the step size
 $\beta^{-1}\in(1,2/L_{\mu \nabla \hat{f}})$.
Hence, from the standard argument of the forward-backward splitting method
(see, e.g., \cite[Example 17.6]{yyy_springer_book11}),
the composition
$\prox_{\beta^{-1}\check{\varphi}}\circ ({\rm Id}- \beta^{-1}\mu \nabla \hat{f})$
is an averaged nonexpansive operator 
with the fixed-point set
$\Fix{\prox_{\beta^{-1}\check{\varphi}}\circ 
({\rm Id}- \beta^{-1}\mu \nabla \hat{f})} 
=\argmin_{x\in\euclidspace} [\mu \hat{f}(x) + \check{\varphi}(x)]
=\argmin_{x\in\euclidspace} [\mu f(x) + \varphi(x)]\neq \varnothing$,
where the nonemptiness is due to the assumption.
The classical Krasnosel'sk\u{i}-Mann iterate
(see Fact \ref{fact:mann})
thus ensures weak convergence of $(x_k)_{k\in\Natural}$
to a fixed point which is a minimizer of $\mu f(x) + \varphi(x)$.
\end{proof}

The existence (and uniqueness) of
minimizer of $\mu f + \varphi$ is ensured
if $\mu\rho> 1-\beta$.
The case in which $f$ is not strongly convex has been
addressed in {\rm \cite{suzuki_icassp24}}.
Note in Theorem \ref{theorem:convergence}.2
that $(x_k)_{k\in\Natural}$ converges 
``strongly (in the norm sense)'' to 
$\hat{x}$ if $\euclidspace$ is finite dimensional
(see Definition \ref{def:convergence}).

\begin{remark}
[On Theorem \ref{theorem:convergence}]
If $\beta\approx 1$ (i.e., $T$ is nearly nonexpansive),
the range of $\mu$ is approximately identical to $(0,2/\kappa)$.
Meanwhile, if $\kappa \approx \rho$,
the range of $\beta$ is nearly $(0,1)$,
meaning that 
the Lipschitz constant $\beta^{-1}$ of $T$ could be 
an arbitrary real number virtually.
However, focusing on the gap with respect to the range of $\mu$,
we observe that
$(1+\beta)/\kappa - (1-\beta)/\rho 
=
K_{\kappa,\rho} (\beta(\kappa+\rho)/(\kappa-\rho) - 1)>0$
with $K_{\kappa,\rho}:= \rho^{-1} - \kappa^{-1}>0$,
where the strict positivity of the gap is verified
by $\beta>(\kappa-\rho)/(\kappa+\rho)$.
This suggests that
the gap of $\mu$ vanishes when 
$\beta^{-1}$ is significantly large so that 
$\beta \approx (\kappa-\rho)/(\kappa+\rho)$.
\end{remark}

\subsection{Primal-Dual Splitting Type Algorithm}
\label{subsec:primal_dual}

We now build a framework which could be applied to a wider class of
problems involving total variation, for instance, 
where a bounded linear operator $L:\euclidspace\rightarrow \mathcal{U}$
will be used in the algorithm.
Here, $\euclidspace$ and $\mathcal{U}$ are real Hilbert spaces equipped
with
the inner products 
$\innerprod{\cdot}{\cdot}_{\euclidspace}$,
$\innerprod{\cdot}{\cdot}_{\mathcal{U}}$,
and the induced norms
$\norm{\cdot}_{\euclidspace}$,
$\norm{\cdot}_{\mathcal{U}}$,
respectively.
We define the operator norm
$\norm{L}:=\sup_{x\in\euclidspace,\norm{x}_{\euclidspace}\leq 1} \norm{Lx}_{\mathcal{U}}$.
The adjoint operator of $L$ is denoted by $L^*:\mathcal{U}\rightarrow \euclidspace$.
With a nonlinear operator $T:\mathcal{U} \rightarrow \mathcal{U}$, 
we consider the following algorithm inspired by the Condat--V\~u
algorithm \cite{condat13,vu13}.\footnote{
An earlier study of an operator-regularization approach using the
primal-dual splitting algorithm can be found in \cite{ono17}.
}

\begin{algorithm}
\label{alg:primal_dual}
 ~\\
{\bf Initialization:} 
$x_0\in \euclidspace$, 
$u_0\in \mathcal{U}$\\
{\bf Requirements:} 
$\sigma>0$,
$\tau>0$,
$\rho >0$

\noindent
$\tilde{u}_{k+1}:= u_k+\sigma L x_k$\\
$u_{k+1}:= \tilde{u}_{k+1} - \sigma T 
\bigg( \bigg(\sigma + \dfrac{\rho}{\norm{L}^2} \bigg)^{-1} \tilde{u}_{k+1} \bigg)$\\
$x_{k+1}:= \bigg({\rm Id} + \dfrac{\tau\rho}{\norm{L}^2} L^* L\bigg)
 x_k- \tau \nabla  f(x_k)- \tau L^*(2 u_{k+1} - u_k)$
\end{algorithm}



\begin{figure*}[t!]
 \psfrag{f}[Bl][Bl][1.0]{$f$ strongly convex}
 \psfrag{g}[Bl][Bl][1.0]{$g$ weakly convex}
 \psfrag{L}[Bl][Bl][1.0]{$L$}
 \psfrag{comp}[Bl][Bl][1.0]{$\circ$}
 \psfrag{smooth}[Bl][Bl][1.0]{smooth}
 \psfrag{+}[Bl][Bl][1.0]{$+$}
 \psfrag{+quad}[Bl][Bl][1.0]{$+$ quadratic}
 \psfrag{-quad}[Bl][Bl][1.0]{$-$ quadratic}
 \psfrag{f-}[Bl][Bl][1.0]{$\hat{f}$ in \eqref{eq:def_f_minus}, convex}
 \psfrag{g+}[Bl][Bl][1.0]{$\check{g}$  in \eqref{eq:def_g_plus}, convex}

 \psfrag{applyPD}[Bl][Bl][1.0]{apply PD}
 \psfrag{recursion}[Bl][Bl][1.0]{\eqref{eq:recursion_x}, \eqref{eq:recursion_u}}
 \psfrag{applyprop}[Bl][Bl][1.0]{apply
 Prop.~\ref{proposition:moreau_decomp_mod} to $\prox_{\sigma \check{g}^*}$
in \eqref{eq:recursion_u}}
 \psfrag{prox}[Bl][Bl][1.0]{\eqref{eq:moreau_decomp_of_gplusconjugate}}
 \psfrag{proxT}[Bl][Bl][1.0]{\eqref{eq:theorem3_desired_identity}}
 \psfrag{replace}[Bl][Bl][1.0]{replace $\sprox_{g/(\sigma+\rho/\norm{L}^2)}$}
 \psfrag{byT}[Bl][Bl][1.0]{by operator $T$}
 \psfrag{substitute}[Bl][Bl][1.0]{substitute}
 \psfrag{algorithm1}[Bl][Bl][1.0]{Algorithm \ref{alg:primal_dual}}
 \psfrag{step 0}[Bl][Bl][1.0]{step 0}
 \psfrag{step 1}[Bl][Bl][1.0]{step 1}
 \psfrag{step 2}[Bl][Bl][1.0]{step 2}
 \psfrag{step 3}[Bl][Bl][1.0]{step 3}
\centering
 \includegraphics[height=4.7cm]{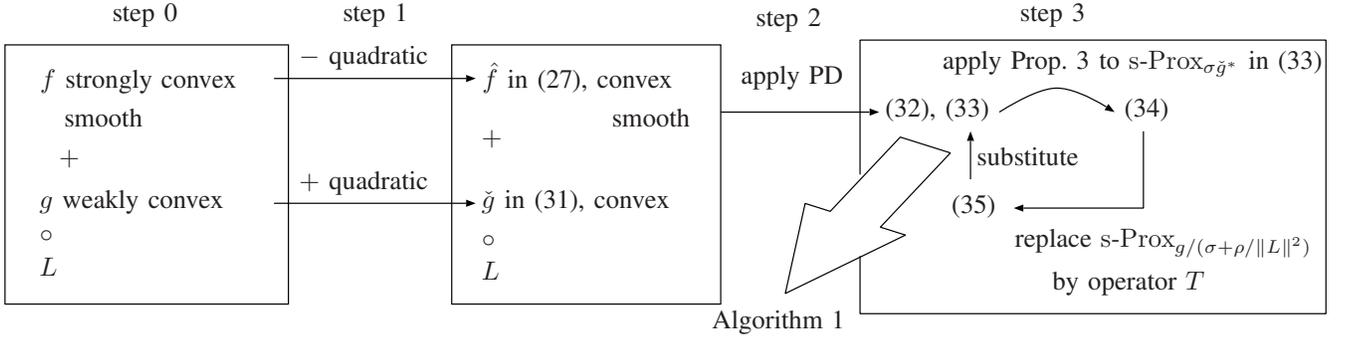}
 \caption{A schematic diagram of the derivation of Algorithm \ref{alg:primal_dual}.}
\vspace{0em}
\label{fig:algorithm1}
\end{figure*}

\begin{theorem}
 \label{theorem:primal_dual}

Let $f\in \gzh$ be a smooth $\rho$-strongly convex function for
 $\rho\in\real_{++}$
so that the function
\begin{equation}
\Gamma_0(\euclidspace)\ni\hat{f}:
x\mapsto f(x) - 
\frac{\rho}{2\norm{L}^2}
\norm{Lx}_{\mathcal{U}}^2
\label{eq:def_f_minus}
\end{equation}
is $\kappa$-smooth for some $\kappa\in\real_{++}$.
Assume that $T:=\nabla \psi:\mathcal{U}\rightarrow \mathcal{U}$
is a MoL-Grad denoiser for some 
$\psi\in\Gamma_0(\euclidspace)$ and
$\beta\in(0,1)$.
Set the step sizes $\sigma,\tau>0$ such that
\begin{itemize}
 \item [(i)]$\sigma\leq 
\dfrac{\rho\beta}{\norm{L}^2(1-\beta)}$,
\item [(ii)] $\tau \Big(\sigma\norm{L}^2 + \dfrac{\kappa}{2} \Big) < 1$.
\end{itemize}
Then, for arbitrary $(x_0,u_0)\in\euclidspace\times \mathcal{U}$,
the sequences $(x_k)_{k\in\Natural} \subset\euclidspace$ and 
$(u_k)_{k\in\Natural} \subset\mathcal{U}$
generated by Algorithm \ref{alg:primal_dual}
converge weakly to solutions $\hat{x}\in\euclidspace$
and  $\hat{u}\in\mathcal{U}$
of the following primal and dual problems, respectively:
\begin{equation}
 \min_{x\in\euclidspace} f(x) + g(Lx),
\label{eq:primal_problem}
\end{equation}
\begin{equation}
 \min_{u\in\mathcal{U}} 
\hat{f}^*(-L^* u) + \check{g}^*(u),
\label{eq:dual_problem}
\end{equation}
provided that such solutions exist.
Here, 
\begin{align}
\hspace*{-2em} g:=&~ \Big(\sigma + \frac{\rho}{\norm{L}^2} \Big)\varphi
= 
\Big(\sigma + \frac{\rho}{\norm{L}^2} \Big)
\Big(\psi^*- \frac{1}{2}\norm{\cdot}_{\mathcal{U}}^2 \Big),
\label{eq:def_g}
\\
\hspace*{-2em} \check{g}:=&~ g + 
\frac{\rho}{2\norm{L}^2}
\norm{\cdot}_{\mathcal{U}}^2 \in
 \Gamma_0(\mathcal{U})
\label{eq:def_g_plus}
\end{align}
 satisfy
 $f+g \circ L = \hat{f} + \check{g} \circ L \in\Gamma_0(\euclidspace)$.
(In this case, $g\circ L$ is the implicit regularizer.)
\end{theorem}

\begin{proof}
First of all,
$\euclidspace\ni x\mapsto (\rho/2)\norm{x}_{\euclidspace}^2 -
 (\rho/(2\norm{L}^2))\norm{Lx}_{\mathcal{U}}^2$ is
a convex function
owing to the positive semidefiniteness of ${\rm Id} - L^* L/\norm{L}^2$.
This, together with the $\rho$-strong convexity of $f$, implies
that $\hat{f} = f - (\rho/(2\norm{L}^2))\norm{L\cdot}_{\mathcal{U}}^2$ is convex.

Since $\varphi$ is $(1-\beta)$-weakly convex by Theorem
 \ref{theorem:weaklyconvex_necsuffcondition}, 
$g=(\sigma + \rho/\norm{L}^2)\varphi$ is
$(\sigma + \rho/\norm{L}^2)(1-\beta)$-weakly convex.
On the other hand, condition (i) of the theorem implies that 
$(\sigma + \rho/\norm{L}^2)(1-\beta)\leq \rho/\norm{L}^2$, 
and hence $\check{g}$ is convex as well.

We show that Algorithm  \ref{alg:primal_dual} is obtained by applying
the Condat--V\~u algorithm to the primal and dual problems,
 \eqref{eq:primal_problem} and \eqref{eq:dual_problem}.
By
$\nabla \hat{f} (x_k)= \nabla f(x_k) - (\rho/\norm{L}^2) L^*L x_k$,
the last recursion of Algorithm \ref{alg:primal_dual} can be rewritten
as
\begin{equation}
 x_{k+1}= x_k- \tau \nabla  \hat{f}(x_k)
 - \tau L^*(2 u_{k+1} - u_k).
\label{eq:recursion_x}
\end{equation}
We shall now show that the first two recursions can be expressed as
\begin{equation}
u_{k+1}= \prox_{\sigma \check{g}^*} (u_k + \sigma L x_k),
\label{eq:recursion_u}
\end{equation}
because the set of recursions 
\eqref{eq:recursion_x} and \eqref{eq:recursion_u}
gives (a no-relaxation version of) the Condat--V\~u iteration (form II) \cite{condat_siam23}.
Indeed, the formulation in \cite{condat_siam23} involves
three functions $f$, $g$, and $h$, each of which corresponds to 
the zero function, $\check{g}$, and $\hat{f}$, respectively, in our case.
Once \eqref{eq:recursion_u} is proved,
the convergence follows immediately by \cite[Theorem 7.1]{condat_siam23}
under condition (ii).

Letting $\phi:=g$ in Proposition \ref{proposition:moreau_decomp_mod} yields
\begin{equation}
 \prox_{\sigma \check{g}^*} = {\rm Id} - \sigma \sprox_{g/(\sigma+
  \rho/\norm{L}^2)} \circ 
\bigg(\sigma+ \frac{\rho}{\norm{L}^2}\bigg)^{-1}
{\rm Id},
\label{eq:moreau_decomp_of_gplusconjugate}
\end{equation}
which, together with
$T = \sprox_{\varphi} =\sprox_{g/(\sigma+ \rho/\norm{L}^2)} $ by
Theorem \ref{theorem:weaklyconvex_necsuffcondition},
verifies
\begin{equation}
 \prox_{\sigma \check{g}^*} = {\rm Id} - \sigma T \circ 
\bigg(\sigma+ \frac{\rho}{\norm{L}^2}\bigg)^{-1}
{\rm Id}.
\label{eq:theorem3_desired_identity}
\end{equation}
The identity in \eqref{eq:theorem3_desired_identity}
justifies the expression in \eqref{eq:recursion_u}, and 
this completes the proof.
\end{proof}

\begin{remark}[Parameter design]
\label{remark:parameter_design}

From \eqref{eq:def_g}, the optimization problem
\eqref{eq:primal_problem} which Algorithm \ref{alg:primal_dual}
actually solves depends on the dual step-size $\sigma$
as well as $\varphi$ which is associated with $T$.
Note here that $\rho$ and $L$ are assumed to be given.
Specifically, the parameter $\sigma$ governs the weight of the regularizer (prior) $\varphi$
relative to $f$, and
a tuning parameter $\delta\in(0,1]$ can be introduced
to control the weight as
$\sigma:= \delta\rho\beta / [\norm{L}^2(1-\beta)]$.
It should be remarked that changing the tuning parameter(s) in $T$ (if any)
may scale $\varphi$ as well as  $\beta$
 which determines the upper bound of $\sigma$ 
(see condition  (i) of Theorem \ref{theorem:primal_dual}).
The primal step size $\tau$ only affects the convergence speed;
it can simply be set to
$\tau:= \gamma(\sigma\norm{L}^2 + \kappa/2)^{-1}$,
where $\gamma\in(0,1)$ (our recommended choice is $\gamma:=0.9$).
\end{remark}




\subsection{A Systematic Way of 
Building Operator-Regularization
	Algorithms using MoL-Grad
with Convergence Guarantee}
\label{subsec:systematic}

We have seen in Section \ref{subsec:primal_dual} that Algorithm
\ref{alg:primal_dual} minimizes $f+g\circ L$ 
(and $\hat{f}^*\circ(-L^*) + \check{g}^*$ at the same time)
under the conditions.
But, {\em how can we derive such an algorithm?}
The derivation process of Algorithm \ref{alg:primal_dual}
is illustrated in Fig.~\ref{fig:algorithm1},
elaborated below in a step-by-step manner.

Suppose that we would like to plug 
a (well-performing) MoL-Grad denoiser $T$
into the primal-dual splitting algorithm.
Since such a $T$ is the proximity operator of a weakly convex function
(Theorem \ref{theorem:weaklyconvex_necsuffcondition}),
we start with the convex objective $f+g\circ L$
with $g$ weakly convex.
First, we define $\hat{f}$ and $\check{g}$ as in 
\eqref{eq:def_f_minus} and \eqref{eq:def_g_plus}, respectively,
which satisfy $f+g\circ L = \hat{f} + \check{g} \circ L$.
Second, the primal-dual algorithm is applied to
the sum $\hat{f} + \check{g} \circ L$ of convex functions 
to obtain the set of recursions
\eqref{eq:recursion_x} and \eqref{eq:recursion_u}.
Third,
Proposition \ref{proposition:moreau_decomp_mod} is applied
to $\prox_{\sigma \check{g}^*}$ to obtain \eqref{eq:moreau_decomp_of_gplusconjugate},
and then $\sprox_{g/(\sigma+  \rho/\norm{L}^2)}$ is replaced by $T$
to obtain \eqref{eq:theorem3_desired_identity}.
Note here that 
$\sprox_{g/(\sigma+  \rho/\norm{L}^2)}$ is a MoL-Grad denoiser
from Theorem \ref{theorem:weaklyconvex_necsuffcondition}
owing to the $(1-\beta)$-weak convexity of 
$g/(\sigma+\rho/\norm{L}^2)(=\varphi)$.
Finally, plugging
\eqref{eq:theorem3_desired_identity}
into \eqref{eq:recursion_u} yields Algorithm \ref{alg:primal_dual}.

The point is that the denoiser can be expressed as
$T=\sprox_{g/(\sigma+ \rho/\norm{L}^2)}$
based on Theorem \ref{theorem:weaklyconvex_necsuffcondition},
which is further linked to the (Moreau's) proximity operator
by Proposition \ref{proposition:moreau_decomp_mod}.
Thanks to this link, Algorithm \ref{alg:primal_dual} has been shown
to  minimize $f+g\circ L$ with the operator-induced function $g$
given in \eqref{eq:def_g}.
In general, 
given an arbitrary operator splitting algorithm,
the proximity operator(s) employed
can be replaced systematically by 
any MoL-Grad denoiser(s)
by the following steps.
\vspace*{.5em}

\noindent{\bf Step 0.}
 Suppose that we are given an operator splitting
	    algorithm in which the proximity operator wants to be
	    replaced by a MoL-Grad denoiser $T$.
The associated optimization problem typically involves
a smooth function $f$ and a prox-friendly convex function $g$ 
(or possibly multiple such functions) 
which is possibly composed with a bounded linear operator.
As opposed to this standard setting, we suppose that $g$ is 
 weakly convex (and possibly other prox-friendly functions could also be so)
 while $f$ is strongly convex so that the whole objective $f+g\circ L$ is convex.

\noindent{\bf Step 1.}
Add a quadratic function to the weakly convex function,
and subtract its corresponding quadratic function from the strongly convex
function in such a way that
(a) the sum remains the same, 
(b) the weakly convex function is convexified, and 
(c) the strongly convex function remains convex.
By doing so, a convex optimization problem is obtained with the same objective.

\noindent{\bf Step 2.}
Apply the operator splitting algorithm to the 
convex optimization problem obtained in step 1.

\noindent{\bf Step 3.}
Apply Proposition \ref{proposition:moreau_decomp_mod}
(or Proposition \ref{proposition:extended_Moreau_decomposition}), and 
replace the proximity operator of $g$ (of some index) by the $T$
(and also replace those of the other functions by other MoL-Grad denoisers).

Through the above steps,
the obtained algorithm involves no proximity operator anymore
in an explicit form.
By the same principle as stated above for 
Algorithm \ref{alg:primal_dual}, however,
it iteratively solves
an optimization problem involving an implicit regularizer(s)
which is induced by the denoiser $T$.



\section{Examples and Relation to Prior Work}\label{sec:example}


We present examples of MoL-Grad denoiser, or equivalently
the proximity operator of weakly convex function
in the sense of Definition \ref{def:sprox}, and we then discuss the relation
to prior work.
In the Euclidean space setting, 
we use the standard notation $\norm{\cdot}_1$ and $\norm{\cdot}_2$
for the $\ell_1$ and $\ell_2$ norms of vector, respectively,
and $(\cdot)^\top$ for matrix/vector transposition.

\subsection{Shrinkage Operators for Scalar}\label{subsec:examples}

Below are a selection of typical shrinkage operators.

\subsubsection{Soft shrinkage}

The soft shrinkage operator 
${\rm soft}_{\lambda}:\real\rightarrow\real$
for the threshold $\lambda \in \real_{++}$
defined by
\begin{equation}
 {\rm soft}(x) :=
\begin{cases}
0 & \mbox{if } \abs{x}\leq\lambda \\
x -\lambda & \mbox{if } x>\lambda\\
x + \lambda & \mbox{if } x< - \lambda
\end{cases}
\end{equation}
is the proximity operator of the convex function 
$\lambda \abs{\cdot}$.
The soft shrinkage operator is a MoL-Grad denoiser for an arbitrary $\beta\in (0,1)$
because it is (firmly) nonexpansive.

\subsubsection{Firm shrinkage}

The firm shrinkage operator 
${\rm firm}_{\lambda_1,\lambda_2}:\real\rightarrow\real$
for the thresholds $\lambda_1\in \real_{++}$
and $\lambda_2 \in (\lambda_1,+\infty)$
defined by
\cite{gao97}
\begin{equation}
\hspace*{-.1em} {\rm firm}_{\lambda_1,\lambda_2}(x) \!:=\!
\begin{cases}
0 & \!\!
\mbox{if } \abs{x}\leq\lambda_1 \\
{\rm sign}(x)\frac{\lambda_2 (\abs{x} - \lambda_1)}{\lambda_2 -
 \lambda_1} & 
\!\!\mbox{if }\lambda_1 < \abs{x} \leq \lambda_2 \\
x & \!\!\mbox{if } \abs{x}>\lambda_2
\end{cases}
\label{eq:firm}
\end{equation}
is the proximity operator $\sprox_\varphi$
of the following $\lambda_1/\lambda_2$-weakly convex function
$\varphi:= \lambda_1 \varphi_{\lambda_2}^{\rm MC}$,
where 
\begin{equation}
\varphi_{\lambda_2}^{\rm MC}(x):= 
\begin{cases}
\abs{x} - \frac{1}{2\lambda_2}x^2 & \mbox{if } \abs{x}\leq \lambda_2 \\
\frac{1}{2}\lambda_2 & \mbox{if } \abs{x}> \lambda_2
\end{cases}
\label{eq:mcp}
\end{equation}
is the minimax concave (MC) penalty \cite{zhang,bayram16,selesnick}.

\subsubsection{Garrote shrinkage}
The garrote shrinkage (nonnegative garrote thresholding) operator 
${\rm garrote}_{\lambda}:\real\rightarrow\real$
for the threshold $\lambda \in \real_{++}$
defined by
\cite{gao98}
\begin{equation}
 {\rm garrote}_{\lambda}(x) :=
\begin{cases}
0 & \mbox{if } \abs{x}\leq\lambda \\
x - \frac{\lambda^2}{x} & \mbox{if } \abs{x}>\lambda
\end{cases}
\end{equation}
is the proximity operator $\sprox_\varphi$
of the following 1/2-weakly convex function:
\begin{align}
\hspace*{-1em} \varphi(x):=  & ~
\frac{1}{4}\left(
\abs{x}\sqrt{x^2 + 4 \lambda^2} - x^2
\right) \nonumber\\
& ~+ 
\lambda^2
\left[
\log (\abs{x} + \sqrt{x^2 + 4\lambda^2}) - \log 2\lambda
\right].
\label{eq:garrote_regularizer}
\end{align}

Note that the hard shrinkage operator 
is {\em not} the proximity operator of any weakly convex function
(in the sense of Definition \ref{def:sprox})
because of its discontinuity.


\subsection{Shrinkage Operators for Vector}\label{subsec:shrinkage_structured}

\begin{definition}
[Moreau enhanced model \cite{selesnick,abe_ip20}]
Given a convex function $f\in\Gamma_0(\euclidspace)$, 
its Moreau enhanced model of index $\lambda\in\real_{++}$
is defined by
\begin{equation}
f_{\lambda}: \euclidspace\rightarrow (-\infty,+\infty]:
x\mapsto f(x) - ~^{\lambda} f(x).
\end{equation}  
\end{definition}

Using this notion, the MC function $\varphi_{\lambda_2}^{\rm MC}$ can be
expressed as the Moreau enhanced model of 
the absolute-value function
$\abs{\cdot}$ of index $\lambda_2$; i.e.,
$\varphi_{\lambda_2}^{\rm MC}(x) =
(\abs{\cdot})_{\lambda_2}$
\cite[Proposition 12]{selesnick}.


\begin{example}[Firm-shrinkage for vector]
\label{example:vector_firm_shrinkage}

Let $\euclidspace$ be an arbitrary Hilbert space equipped with a norm
$\norm{\cdot}$.
For positive constants $\lambda_1\in \real_{++}$ and $\lambda_2 \in (\lambda_1,+\infty)$,
we define the vector firm-shrinkage operator:
\begin{equation}
 T^{\rm v\text{-}firm}_{\lambda_1,\lambda_2}:\euclidspace \rightarrow \euclidspace:
(0\neq )x\mapsto \frac{x}{\norm{x}}
{\rm firm}_{\lambda_1,\lambda_2}(\norm{x})
\label{eq:Tlambda1lambda2}
\end{equation}
with $T^{\rm v\text{-}firm}_{\lambda_1,\lambda_2}(0):=0$.
Then, $T^{\rm v\text{-}firm}_{\lambda_1,\lambda_2}$ is the proximity operator
of the ($\lambda_1/\lambda_2$)-weakly convex function $\lambda_1 (\norm{\cdot})_{\lambda_2}$,
where
$ (\norm{\cdot})_{\lambda_2}
(:=\norm{\cdot} - ~^{\lambda_2}\!\norm{\cdot})
= \varphi_{\lambda_2}^{\rm MC}
\circ \norm{\cdot}$.
\end{example}

\begin{example}[Group firm-shrinkage]
\label{example:group_firm_shrinkage}

Let $\euclidspace_1$, $\euclidspace_2$, $\cdots$, $\euclidspace_G$
be real Hilbert spaces
equipped with the norms $\norm{\cdot}_{\euclidspace_1}$, $\norm{\cdot}_{\euclidspace_2}$, $\cdots$, $\norm{\cdot}_{\euclidspace_G}$,
respectively, and define the product space
$\euclidspace:=   \euclidspace_1 \times \euclidspace_2
 \times\cdots\times \euclidspace_G$
equipped with the norm
$\norm{x}_\euclidspace:=\sqrt{
\sum_{i=1}^{G}\norm{x_i}_{\euclidspace_i}^2}$,
$x:= (x_i)_{i=1}^G\in\euclidspace$.
Given constants $\lambda_1\in \real_{++}$ and $\lambda_2 \in (\lambda_1,+\infty)$,
we define the group firm-shrinkage operator:
\begin{equation}
 T^{\rm g\text{-}firm}_{\lambda_1,\lambda_2}:\euclidspace\rightarrow \euclidspace:
 (x_i)_{i=1}^G \mapsto
\left(
T^{\rm v\text{-}firm}_{\lambda_1,\lambda_2}(x_i)
\right)_{i=1}^G,
\end{equation}
where $T^{\rm v\text{-}firm}_{\lambda_1,\lambda_2}$ is defined by 
\eqref{eq:Tlambda1lambda2} in each space.
Then, $T^{\rm g\text{-}firm}_{\lambda_1,\lambda_2}$ 
is the proximity operator of the $(\lambda_1/\lambda_2)$-weakly convex function
$\lambda_1 (\norm{\cdot}_{{\rm g},1})_{\lambda_2}$, where
$\norm{\cdot}_{{\rm g},1}: 
\euclidspace\rightarrow \real_+:
 (x_i)_{i=1}^G \mapsto \sum_{i=1}^{G} \norm{x_i}_{\euclidspace_i}$.

\end{example}

The proofs of Examples \ref{example:vector_firm_shrinkage} and
\ref{example:group_firm_shrinkage} are given in
 Appendices 
\ref{sec:proof_vector_firm} and \ref{sec:proof_group_firm},
respectively.

\begin{example}[Neural Network with Tied Weights]
\label{example:tied_weights}

Let $\nabla\psi:\real^m \rightarrow \real^m$ be 
a firmly-nonexpansive activation operator
(see Section \ref{subsec:nonexpansive}).
In fact, many of the known activation functions are firmly nonexpansive
{\rm \cite{combettes_dnn20}}, such as ReLU, sigmoid, and softmax.
For instance, ReLU 
$\real\rightarrow \real:
x\mapsto
 \begin{cases}
  x & \mbox{ if } x\in \real_+\\
0 & \mbox{ otherwise}
 \end{cases}$
can be expressed as
the proximity operator $\prox_{\varphi}$
of the indicator function $\varphi:=\iota_{\real_+}:
\real\rightarrow \{0,+\infty\}:
x\mapsto 
 \begin{cases}
  0 & \mbox{ if } x\in \real_+\\
+\infty & \mbox{ otherwise,}
 \end{cases}
$
where $\prox_{\varphi}$ is firmly nonexpansive owing to the convexity of 
$\iota_{\real_+}$ (see Section \ref{subsec:prox}).
Let $\signal{W}\in\real^{m\times n}$
be a nonzero weight matrix.
Then, a weight-tied neural network {\rm \cite{vincent10}}
$T:= 
\signal{W}^\top \circ \nabla\psi\circ \signal{W}
= \nabla (\psi\circ \signal{W})$
can be expressed as
$T= \sprox_{\varphi}$ with
$\varphi:= (\psi\circ \signal{W})^* - (1/2)\norm{\cdot}_2^2$.
Here, if $\kappa:= \|\signal{W}^\top\signal{W}\|>1$, 
$T$ is a MoL-Grad denoiser, and $\varphi$ is $(1-\kappa^{-1})$-weakly
 convex.
If $\kappa\in (0,1]$, in contrast,
$T$ is Moreau's proximity operator of
the convex function $\varphi$.
More practical network architectures
are discussed in {\rm \cite{shimizu_sip24}},
where multi-layered neural networks with tied nonnegative weights
are shown to be MoL-Grad denoisers.
\end{example}




\begin{example}[Debiased OSCAR \cite{suzuki_icassp24}]
\label{example:doscar}

Given constants $\lambda_1,\lambda_2\in \real_{++}$,
the octagonal shrinkage and clustering algorithm for regression (OSCAR)
{\rm \cite{bondell08,zhong12,zeng15}}
is defined by
\begin{equation*}
 \Omega_{\lambda_1,\lambda_2}^{{\rm OSCAR}}:\real^n\rightarrow \real:
\signal{x}\mapsto \lambda_1\norm{\signal{x}}_1 + \lambda_2
\sum_{i<j} \max \{\abs{x_i},\abs{x_j}\}.
\end{equation*}
For $\omega,\eta\in (1,+\infty)$, the operator
\begin{equation}
T_{\lambda_1,\lambda_2,\omega,\eta}^{{\rm DOSCAR}}
:= \omega \prox_{\Omega_{\lambda_1,\lambda_2}^{{\rm OSCAR}}}
 - (\omega -1)\prox_{\eta \Omega_{\lambda_1,\lambda_2}^{{\rm OSCAR}}}
\end{equation}
is called the debiased OSCAR (DOSCAR)
{\rm \cite{suzuki_icassp24}}.
The operator $T_{\lambda_1,\lambda_2,\omega,\eta}^{{\rm DOSCAR}}$ is 
the proximity operator of a certain weakly convex function.
\end{example}

Finally, we mention that a principled way of deriving a continuous
relaxation
of a given discontinuous shrinkage operator is discussed 
under the MoL-Grad framework in \cite{yukawa_erowl24},
giving another example of the MoL-Grad denoiser.



\subsection{Prior Art I: Plug-and-Play Methods}
\label{subsec:relation_to_prior_work}

\noindent {\bf Convergence analysis of plug-and-play (PnP) method:}
There exist a significant amount of works concerning 
the theoretical properties of the PnP method
(see \cite{ahmad20,mukherjee23,kamilov23} for extensive lists of references).
Among many others, \cite{ryu19} has proved the convergence of the PnP
forward-backward splitting (PnP-FBS) method and
the PnP Douglas-Rachford splitting (PnP-DRS) method,
where every single iterate is assumed contraction by restricting
the Lipschitz constant of ${\rm Id} - T$
with the so-named real spectral normalization (realSN).
Here, the Lipschitz constant must be smaller, for instance, than
$\varepsilon:= 2\rho/(\kappa-\rho)$
in the case of the PnP-FBS method (see Theorem \ref{theorem:convergence}
for definition of $\rho$ and $\kappa$).
In this framework, the limit point of the iterate is the unique fixed
point of the contraction mapping (Banach-Picard).

In \cite{sun19}, the class of averaged-nonexpansive denoiser has been studied,
encompassing those denoisers which cannot be expressed in terms of the
proximity operator.
In this framework, the limit point is characterized as a
fixed point of the iterate.
In those works, the implicit regularizer (if exists) has not been
given in an explicit form.
There are several previous works considering the implicit regularizer:
the so-called gradient-step denoiser and the regularization by denoising
approach.
Those two approaches will be discussed below.

\begin{table}
\caption{Comparisons of the existing denoisers for the plug-and-play
 method in terms of
(i) the condition on the Lipschitz constant L$(\cdot)$,
(ii) smoothness of the implicit regularizer $\varphi$, and
(iii) the range of the weak-convexity constant of $\varphi$}
\label{table:pnp}

\centering
\begin{tabular}{c|c|c|c}
method & Lipschitz condition & smoothness & range of ($1-\beta$) \\
 \hline 
realSN \cite{ryu19} & L$({\rm Id} - T) < \varepsilon$
& - & - \\ \hline
averaged \cite{sun19} & L$(T)\leq 1$ & - & - \\ \hline
GS \cite{hurault22} & L$({\rm Id} -T) < 1$ & smoooth &
$(0,1/2)$
\\ \hline
RED \cite{romano17} & L$(T)\leq 1$ & - & - \\ \hline
MMO \cite{pesquet21} & L$(2T - {\rm Id})\leq 1$ & - & - \\ \hline
MoL-Grad & {\bf free} & {\bf free} & 
$(0,1)$
\\ \hline
\end{tabular}

\end{table}

\noindent {\bf Gradient-step (GS) denoiser:}
The GS denoiser has been studied in the literature
\cite{cohen21,hurault_iclr22,hurault22}.
Among those studies, the most relevant results 
would be the one in \cite[Proposition 3.1]{hurault22},
where $T:={\rm Id} - \nabla g$ is assumed to be contraction
($\kappa$-Lipschitz continuous for $\kappa<1$).\footnote{
The arguments of \cite[Proposition 3.1]{hurault22} rely essentially on
\cite{gribonval11}.
}
Because of this restrictive assumption,
the result in \cite[Proposition 3.1]{hurault22} has several limitations:
(i) $T$ possesses a positive definite Jacobian matrix at every point, and
(ii) $\varphi$ is smooth and $\eta_{\rm GS}$-weakly convex for $\eta_{\rm GS}\in(0,1/2)$
($\psi:=(1/2)\norm{\cdot}^2 -g$ is strongly convex and smooth at the same time).
In particular, 
the smoothness of $\varphi$ could be a strict limitation
(see Remark \ref{remark:mol_grad_smoothness} below).
The implicit regularizer given in \cite[Proposition
3.1]{hurault22} can be reproduced 
as a special case of Theorem \ref{theorem:weaklyconvex_necsuffcondition}
straightforwardly.

In \cite{xu20}, the convergence of the PnP method with
the minimum mean squared error (MMSE) denoiser has been studied,
where the MMSE denoiser is expressed in a form of GS denoiser
through {\em Tweedie's formula}.
The implicit regularizer is thus given in the same way as in 
\cite[Proposition 3.1]{hurault22}, and 
it must be smooth accordingly.

\noindent {\bf Regularization by denoising (RED):}
The RED approach \cite{romano17} relates specific
algorithms to a variational optimization problem.
We emphasize that it is {\em not} such an approach that relates a denoiser
itself directly to an implicit regularizer.
In this approach, 
the denoiser $T:\real^n\rightarrow \real^n$ is assumed to possess
(i) nonexpansiveness,
(ii) symmetric Jacobian $J T(\signal{x})$ at every point $\signal{x}\in\real^n$, and
(iii) the so-called local
homogeneity: $T(c\signal{x}) = c T(\signal{x})$ for every $c\in [1-\epsilon,1+\epsilon]$
for a sufficiently small $\epsilon>0$.
Here, (iii) is a key assumption to ensure
$T(\signal{x}) = [J T(\signal{x})] \signal{x}$, from which
it can be shown that 
$\nabla h_{\rm RED}(\signal{x}) = \signal{x} - T(\signal{x})$ for
$h_{\rm RED}(\signal{x}):= (1/2)\signal{x}^{\top} (\signal{x} - T(\signal{x}))$.
We stress here that the relation between
the denoiser $T$ and the regularizer $h_{\rm RED}$
is essentially different from our case.
It has been pointed out that the RED algorithm does not actually
minimize the variational objective \cite{reehorst19}, and another
approach called {\em consensus equilibrium} has also been proposed
as an optimization-free generalization \cite{buzzard18}.

\noindent {\bf Characterization of the limit point for plug-and-play method:}
In {\rm \cite[Proposition 3.1]{pesquet21}},
a characterization of the limit point 
of the vector sequence generated by the forward-backward splitting
algorithm using the plug-and-play method has been
presented through variational inequalities
for a general family of maximally monotone operator.
Their analysis covers a general case rather than focusing solely on
optimization.
However,
the operator is supposed to be firmly nonexpansive in 
{\rm \cite[Proposition 3.1]{pesquet21}}, whereas
it is not necessarily nonexpansive in our case.
In addition, we explicitly showed in Section \ref{subsec:primal_dual}
that, when the plug-and-play is applied to the primal-dual splitting algorithm,
it is still possible to ensure the convergence to an optimal point
with a slight modification.
\begin{remark}[Novelty of MoL-Grad denoiser]
\label{remark:mol_grad_smoothness}
Table \ref{table:pnp} shows comparisons of the proposed and existing approaches.
To the best of the authors' knowledge, MoL-Grad is the first denoiser
that requires no control of its Lipschitz constant to guarantee the
convergence of the operator-regularization (PnP) algorithms.
In addition, the implicit regularizer 
$\varphi(=\psi^* - (1/2)\norm{\cdot}^2)$ 
induced by MoL-Grad denoisers
is not limited to smooth ones.
This means that the class of MoL-Grad denoiser is sufficiently large
to contain those denoisers inducing nonsmooth regularizers.
Nonsmoothness of regularizers has actually played a crucial role in sparse
 modeling.
For instance, the firm/garrote shrinkage presented in 
Section \ref{subsec:examples} induces the nonsmooth regularizer
given in \eqref{eq:mcp} or \eqref{eq:garrote_regularizer}.
For the weight-tied neural networks discussed in Example
\ref{example:tied_weights}, moreover,
adoption of a saturating activation function in the network such as ReLU
makes the associated regularizer nonsmooth,
because
$\nabla^2 \psi(\signal{y})= J T(\signal{y}) = \signal{O}$ at some $\signal{y}\in\real^m$.

We finally mention that, for $\varsigma\in\real_{++}$,
(i) $\varphi$ is $\varsigma$-smooth 
$\Rightarrow$ 
(ii) $\varphi + (1/2)\norm{\cdot}^2 (= \psi^*)$ is 
$(\varsigma+1)$-smooth
$\Leftrightarrow$ 
(iii) $\psi$ is $1/(\varsigma+1)$-strongly convex 
$\Rightarrow$ 
(iv) $T=\nabla \psi$ is 
$1/(\varsigma+1)$-strongly monotone
(see Section \ref{subsec:convex} for the equivalence (ii)
 $\Leftrightarrow$  (iii)).
This implies that the GS denoiser is restricted to strongly monotone operators,
while the MoL-Grad denoiser is free from such a restriction.
\end{remark}

\subsection{Prior Art II: Other Related Works}
\label{subsec:relation_to_prior_work2}

\noindent {\bf Characterization of implicit regularizer:}
In {\rm \cite{combettes_dnn20}},
the softmax activation function is shown to be Moreau's proximity
operator of a certain convex function
by using {\rm \cite[Corollary 24.5]{combettes}},
which is related to, but is different from,
Theorem \ref{theorem:weaklyconvex_necsuffcondition}(b).
Specifically, 
{\rm \cite[Corollary 24.5]{combettes}} states that,
given $\varphi\in\Gamma_0(\euclidspace)$ and
$f\in\Gamma_0(\euclidspace)$
such that $\varphi = f - (1/2)\norm{\cdot}^2$, it holds that
 $\sprox_{\varphi}= \nabla f^*$.
Here, $f$ is implicitly assumed to be $1$-strongly convex.
Meanwhile, by changing
$\psi^*$ in Theorem \ref{theorem:weaklyconvex_necsuffcondition}(b)
by $f$,
it can essentially be stated as follows:
given
$\varphi + ((1-\beta)/2)\norm{\cdot}^2\in\Gamma_0(\euclidspace)$
and 
$f - (\beta/2)\norm{\cdot}^2\in\Gamma_0(\euclidspace)$
with $\varphi = f - (1/2)\norm{\cdot}^2$,
it holds that $\sprox_\varphi = \nabla f^*$.
Here, $\varphi$ is $(1-\beta)$-weakly convex 
(while $f$ is $\beta$-strongly convex) in contrast to the former case.


\noindent {\bf Proximity operator of nonconvex function:}
In {\rm \cite{bauschke21}},
it has been shown  that,
given a proper lower semicontinuous function,
it is weakly convex if and only if 
its proximity operator is cocoercive.
The arguments therein are based on
the notion of {\it abstract subgradient},
while our arguments are based solely on the standard subgradient 
(see \eqref{eq:def_subdifferential}) adopted in convex analysis.
In {\rm \cite[Proposition 2]{gribonval20}}, 
given a mapping $T:\euclidspace\rightarrow \euclidspace$ and $L>0$,
the equivalence of the following two statements has been established:
(a) $T(x)\in \mathbf{Prox}_\varphi(x)$, $\forall x\in \euclidspace$,
for some $\varphi:\euclidspace\rightarrow (-\infty,+\infty]$ such that
$\varphi+ (1/2)[1-1/L] \norm{\cdot}^2\in\Gamma_0(\euclidspace)$
(see Appendix \ref{subsec:proof_prop1} 
for the definition of set-valued operator $\mathbf{Prox}_\varphi$),
and
(b) $T$ is $L$-Lipschitz continuous and
$T(x)\in \partial \psi(x)$, $\forall x\in \euclidspace$,
for some $\psi\in\Gamma_0(\euclidspace)$.
In {\rm \cite{gribonval20}},
the proximity operator is defined as a set-valued operator (or its selection),
and the optimization aspect has not been discussed.

In contrast to the studies in \cite{bauschke21,gribonval20},
our s-prox operator is defined 
as a ``unique'' minimizer of the penalized function (see Definition \ref{def:sprox}).
This is because our primal focus is on the explainability
 perspective of optimization algorithms, and because
for this reason our denoiser is continuous
(see Proposition \ref{proposition:single_valuedness}).
Theorem \ref{theorem:weaklyconvex_necsuffcondition} presented in
Section \ref{subsec:wc_and_coco} is a refinement of 
{\rm \cite[Proposition 2]{gribonval20}}.
Specifically, Theorem \ref{theorem:weaklyconvex_necsuffcondition}
explicitly shows that $\mathbf{Prox}_\varphi$ is actually single-valued,
and it discloses the exact relation between $\varphi$ and $\psi$, based
on which the characterization of implicit regularizers is given
for the operator-regularization approaches
in Section \ref{sec:applications}.

\noindent {\bf Proximity operator of weakly convex function in a linear
inverse problem:}
In {\rm \cite{bayram16}}, the proximity operator 
(in the sense of Definition \ref{def:sprox} essentially)
of weakly convex functions  has been studied in a linear inverse problem,
and the convergence of the forward-backward splitting algorithm with 
the s-prox operator has been analyzed.
While the main results therein are related to Theorem \ref{theorem:convergence},
the present study has a wider scope.
Specifically, it includes 
Theorem \ref{theorem:weaklyconvex_necsuffcondition},
Propositions \ref{proposition:extended_Moreau_decomposition}
and \ref{proposition:moreau_decomp_mod}, and
Theorem \ref{theorem:primal_dual} as well as the discussions in Section
\ref{subsec:systematic} to make the idea be applicable to other 
operator splitting algorithms.
Moreover, the study in \cite{bayram16} 
considers the case in which
(i) $\euclidspace$ is the Euclidean space,
and (ii) the range of $\varphi$ is $\real$.
whereas the present study concerns the case in which
(i) $\euclidspace$ is a general Hilbert space,
and 
(ii) $\varphi$ is allowed to take the value $+\infty$.

\noindent {\bf Other works:}
There are some other related works.
In {\rm \cite{chartrand14}}, a class of ``separable'' shrinkage operators
(i.e., one-dimensional shrinkage operators essentially)
and their induced penalty functions have been studied.
In {\rm \cite{bayram15}}, a study on one-dimensional monotone operators
has been presented, showing that
``a non-decreasing non-constant function with at most a countable
number of discontinuities'' can be expressed as 
a ``selection'' of the (set-valued) proximity operator
of a weakly convex function.
Although those studies are related to Theorem \ref{theorem:weaklyconvex_necsuffcondition},
the scope of the present study is different in the sense that
it concerns the operator-regularization approaches explicitly
based on the s-prox operator of nonseparable weakly-convex function
defined on a Hilbert space.

\section{Simulations}\label{sec:numerical}

Simulations are conducted to show (i) how the theory works in practice
and (ii) how the generalized proximity operator competes with 
Moreau's one.
For this purpose, the firm shrinkage operator will be adopted
as the MoL-Grad denoiser $T$ in Algorithm \ref{alg:primal_dual}.

\subsection{Firm Shrinkage Plugged into Algorithm \ref{alg:primal_dual}}
\label{subsec:algorithm1_with_firm}

We consider the linear system model
$\signal{y}:=\signal{A}\signal{x}_{\diamond}
+\signal{\varepsilon}\in\real^m$,
where 
$\signal{x}_{\diamond}\in \euclidspace := \real^n$
is the $n$-dimensional Euclidean vector to estimate,
$\signal{A}\in\real^{m\times n}$ is the sensing matrix,
and $\signal{\varepsilon}\in\real^m$ is the noise vector.
We consider the overdetermined case
$\lambda_{\min}(\signal{A}^\top\signal{A})>0$,
where $\lambda_{\min}(\cdot)$ denotes the smallest eigenvalue.
(See \cite{suzuki_icassp24} for the case of
$\lambda_{\min}(\signal{A}^\top\signal{A})=0$.)
We suppose that $\signal{x}_{\diamond}$ is piecewise constant
so that the difference vector
$\signal{D}\signal{x}=[x_1-x_2,x_2-x_3,\cdots,x_{n-1}-x_n]^\top$ is
sparse,
where $\signal{D}:=[\signal{I} ~ \signal{0}] -
[\signal{0} ~ \signal{I}]\in \real^{(n-1)\times n}$
with the identity matrix
$\signal{I}$
and the zero vector
$\signal{0}$ of length $n-1$.
We therefore let 
$L : \real^n \rightarrow \real^{n-1}=:\mathcal{U}:
\signal{x} \mapsto \signal{D}\signal{x}$.
The data fidelity term is set to 
$f(\signal{x}):= (1/2)\norm{\signal{A}\signal{x} -
\signal{y}}^2_2$.
In this case, 
since the Hessian matrices of $f$ and $\hat{f}$ are given by 
$\nabla^2 f(\signal{x}) = \signal{A}^\top\signal{A}$ 
and 
$\nabla^2 \hat{f}(\signal{x}) = \signal{A}^\top\signal{A} - 
(\rho/ \norm{\signal{D}}^2)\signal{D}^\top \signal{D}$, respectively,
at every $\signal{x}\in\real^n$,
$f$ is $\rho$-strongly convex
for $\rho:= \lambda_{\min}(\signal{A}^\top\signal{A})$,
and $\hat{f}$ is $\kappa$-smooth for
$\kappa := \|\signal{A}^\top\signal{A} - 
(\rho/ \norm{\signal{D}}^2)\signal{D}^\top \signal{D}\|$.

Recall that the firm shrinkage given in \eqref{eq:firm} has two
parameters
$\lambda_1$ and $\lambda_2$, and its corresponding weakly convex
function is $\varphi:=\lambda_1 \varphi_{\lambda_2}^{\rm MC}$
with $\varphi +
(\lambda_1/(2\lambda_2))\norm{\cdot}^2\in\Gamma_0(\real^n)$,
which means that $\beta= 1-\lambda_1/\lambda_2$.
Hence, in view of Theorem \ref{theorem:primal_dual},
under an appropriate setting,
the sequence $(\signal{x}_k)_{k\in\Natural}$ generated by Algorithm
\ref{alg:primal_dual} converges to a solution
of the following minimization problem:
\begin{equation}
\min_{\signal{x}\in\real^n} f(\signal{x}) +
[\sigma + 
  \lambda_{\min}(\signal{A}^\top\signal{A})/
  \norm{\signal{D}}^2
] \lambda_1 \varphi_{\lambda_2}^{\rm MC} (\signal{D}\signal{x}).
\label{eq:firm_cost_function}
\end{equation}
From \eqref{eq:firm_cost_function}, 
the impact of the regularizer $\varphi_{\lambda_2}^{\rm MC}$ scales
with $\lambda_1$ as well as the dual step size $\sigma$.

As mentioned in Remark \ref{remark:parameter_design},
the $\sigma$ satisfying 
condition (i) in Theorem \ref{theorem:primal_dual}
can be parametrized
by $\delta\in(0,1]$ as 
$\sigma:= 
\delta 
(
  \lambda_{\min}(\signal{A}^\top\signal{A})/
  \norm{\signal{D}}^2
) 
[(\lambda_2/\lambda_1) - 1]>0$,
making  \eqref{eq:firm_cost_function} into
\begin{equation}
\min_{\signal{x}\in\real^n} f(\signal{x}) +
\mu^{-1} \varphi_{\lambda_2}^{\rm MC} (\signal{D}\signal{x}),
\label{eq:firm_cost_function2}
\end{equation}
where
$\mu:= \norm{\signal{D}}^2/
 [ \lambda_{\min}(\signal{A}^\top\signal{A})((1-\delta) \lambda_1 +
 \delta \lambda_2)] >0$.
Here, 
$\norm{\signal{D}}^2$ and
$\lambda_{\min}(\signal{A}^\top\signal{A})$
serve as a sort of normalization,
while $(1-\delta) \lambda_1 + \delta \lambda_2\in (\lambda_1,\lambda_2]$
weights the regularizer $\varphi_{\lambda_2}^{\rm MC}$
relative to $f$; 
increasing $\delta$ shifts
the weight from $\lambda_1$ to $\lambda_2$.

Letting $\delta:=1$,
or equivalently setting $\sigma$ to its upper bound,
reduces  \eqref{eq:firm_cost_function} to
\begin{equation}
\min_{\signal{x}\in\real^n} f(\signal{x}) +
 \frac{\lambda_2 \lambda_{\min}(\signal{A}^\top\signal{A})}
{\norm{\signal{D}}^2}
  \varphi_{\lambda_2}^{\rm MC} (\signal{D}\signal{x}),
\label{eq:firm_cost_function_upper_bound}
\end{equation}
which is independent of $\lambda_1$.
Note here that the dependence on $\lambda_1$ is eliminated by 
condition (i) to guarantee convexity of the overall cost $f+g\circ L$.
As such, the performance of Algorithm \ref{alg:primal_dual} in this case
is governed solely by $\lambda_2$, given the matrices $\signal{A}$
and $\signal{D}$.
This might be advantageous in practice, because $\lambda_2$ is the only
parameter to tune, unless one cares the speed of convergence.


\subsection{Experimental Verification of Theorem \ref{theorem:primal_dual}}
\label{subsec:verification_of_theorem}

Since $\varphi_{\lambda_2}^{\rm MC}
= \norm{\cdot}_1 - ~^{\lambda_2}(\norm{\cdot}_1)$,
\eqref{eq:firm_cost_function2}
can be rewritten as follows:
\begin{align}
\min_{\signal{x}\in\real^n} 
\underbrace{
   \mu f(\signal{x}) 
   -
     ~^{\lambda_2}(\norm{\cdot}_1)
    (\signal{D} \signal{x})
}_{=\tilde{f}(\signal{x})}
+
\norm{
\signal{D} \signal{x}
}_1.
\label{eq:firm_cost_function_with_mu}
\end{align}
As $\tilde{f}:=   \mu f
   -  ~^{\lambda_2}(\norm{\cdot}_1)
    \circ \signal{D}$ is a smooth convex function 
with its gradient available in a closed form,
the minimization problem in \eqref{eq:firm_cost_function_with_mu}
can be solved by the operator splitting algorithm such as 
the primal-dual splitting method with convergence guarantee.\footnote{
This is possible in this specific example, because the function
$\varphi$ corresponding to the firm shrinkage is available.
In general, $\varphi$ could be unavailable, and thus this kind of
approach cannot be used in such a case.
}
To demonstrate the validity of Theorem \ref{theorem:primal_dual} by simulation,
we apply the Condat--V\~u algorithm (form II) to
\eqref{eq:firm_cost_function_with_mu} to see whether Algorithm
\ref{alg:primal_dual} converges to the same point.
We show, more specifically, that
the discrepancy is vanishing
between 
 $(\signal{x}_k,\signal{u}_k)$ of Algorithm \ref{alg:primal_dual} 
and 
 $(\tilde{\signal{x}}_k,\tilde{\signal{u}}_k)$ of the Condat--V\~u
 algorithm
 applied to \eqref{eq:firm_cost_function_with_mu},
where ``discrepancy'' is quantified as follows:
\begin{equation}
\frac{ \norm{\signal{x}_k - \tilde{\signal{x}}_k}_2^2
+  \norm{\signal{u}_k - \tilde{\signal{u}}_k}_2^2}
{ \norm{\signal{x}_k }_2^2
+  \norm{\signal{u}_k}_2^2
}.
\label{eq:discrepancy}
\end{equation}

\begin{figure}[t!]
 \psfrag{xk}[Bl][Bl][0.9]{$\signal{x}_{\diamond}$}
 \psfrag{k}[Bc][Bc][0.9]{index}
\centering
 \includegraphics[width=7cm]{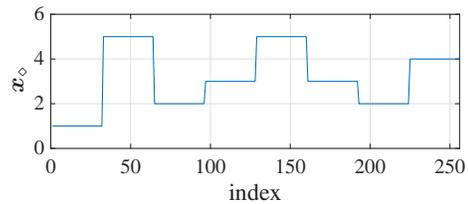}
 \caption{ The vector $\signal{x}_{\diamond}$ used in the simulation.}
\label{fig:xtrue}
\end{figure}

\begin{figure}[t!]
 \psfrag{x}[Bl][Bl][0.9]{$\signal{x}$}
 \psfrag{u}[Bl][Bl][0.9]{$\signal{u}$}
 \psfrag{x, u}[Bl][Bl][0.9]{$\signal{x}, \signal{u}$}
\centering
 \includegraphics[width=8cm]{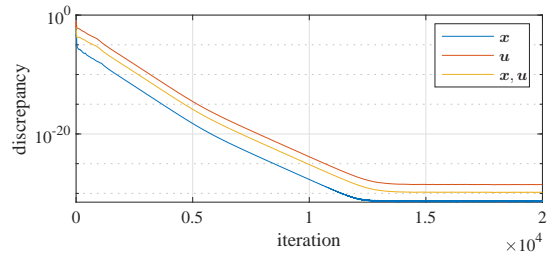}
 \caption{Agreement of the learning curves in line with Theorem \ref{theorem:primal_dual}.}
\label{fig:confirm}
\end{figure}

We let $n:=256$, $m=1024$, and $\signal{x}_{\diamond}$ depicted in
Fig.~\ref{fig:xtrue} is used.
For Algorithm \ref{alg:primal_dual},
we set $\delta:=1$, $\tau:= 0.9(\sigma\norm{\signal{D}}^2 + \kappa/2)^{-1}$,
$\lambda_1:=2.5$, and $\lambda_2:=5$.
We repeat that $\lambda_1$ does not change the solution, see
\eqref{eq:firm_cost_function_upper_bound}.
For the approach based on \eqref{eq:firm_cost_function_with_mu},
we set $\sigma:=0.2$, and all other parameters are the same as
for Algorithm \ref{alg:primal_dual}.
Figure \ref{fig:confirm} plots the discrepancy
given in \eqref{eq:discrepancy} (which is labeled as `$\signal{x}, \signal{u}$');
the discrepancy in terms of the primal variables
and that in terms of the dual variables are also plotted, 
labeled as `$\signal{x}$' and `$\signal{u}$', respectively.
It can be seen that the discrepancy vanishes as time goes by.

\begin{figure}[t!]
 \psfrag{x}[Bl][Bl][0.9]{$\signal{x}$}
 \psfrag{u}[Bl][Bl][0.9]{$\signal{u}$}
 \psfrag{x, u}[Bl][Bl][0.9]{$\signal{x}, \signal{u}$}
\centering

 \includegraphics[width=8cm]{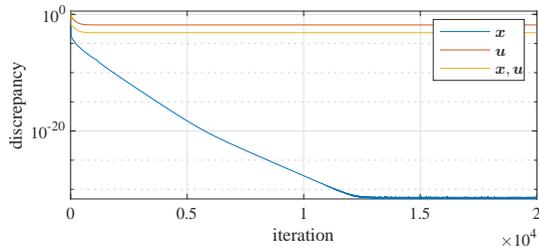}

 \caption{Disagreement of the learning curves for the straightforward approach.}
\label{fig:disagreement}
\end{figure}

For a comparison, we also test a straightforward approach
(which has no theoretical guarantee):
apply the Condat-V\~u algorithm (form II)
directly to the nonconvex formulation
\eqref{eq:firm_cost_function2},
which is essentially the same as \eqref{eq:firm_cost_function_with_mu}.
Here, $\prox_{\sigma g^*}$ appearing in the algorithm
is replaced\footnote{
This replacement is heuristic, and there is no theoretical support
for that.
} by ${\rm Id} - \sigma 
\sprox_{\varphi_{\lambda_2}^{\rm MC}/(\mu \sigma)} \circ(\sigma^{-1}{\rm Id})$,
where
$\sprox_{\varphi_{\lambda_2}^{\rm MC}/(\mu \sigma)}$ is the proximity
operator in the sense of Definition \ref{def:sprox},
which coincides with
the firm shrinkage with $\lambda_1:=1/(\mu \sigma)$.
The $\sigma$ here needs to be chosen in such a way that $\lambda_1 < \lambda_2$.
The results are plotted in Fig.~\ref{fig:disagreement}.
It can be seen that the discrepancy (`$\signal{x}, \signal{u}$') does not vanish.
Although  the discrepancy in terms of `$\signal{x}$' vanishes,
there is no theoretical guarantee for that,
to the best of authors' knowledge.

\subsection{Performance Comparisons}
\label{subsec:performance_comp}

We compare the performance of
the approach described in
Section \ref{subsec:algorithm1_with_firm}
to that of the $\ell_1$-based total variation approach
under the same setting basically as in 
Section \ref{subsec:verification_of_theorem}.
For the former approach, 
since $\lambda_2$ is the only parameter that affects
the solution, we change $\lambda_2$ to see how the performance changes
accordingly.
The cost function of the latter approach is given by
$\mu f + \norm{\cdot}_1\circ \signal{D}$, and we change $\mu$ for this
approach.
The evaluation metric is system mismatch defined by
$\norm{\signal{x}_{\diamond} - \signal{x}_k}_2^2/\norm{\signal{x}_{\diamond}}_2^2$.

\begin{figure}[t!]
 \psfrag{mu}[Bl][Bl][0.9]{$\mu$ ($\delta$)}
 \psfrag{L1}[Bl][Bl][0.9]{$\ell_1$}

\centering

 \includegraphics[width=8cm]{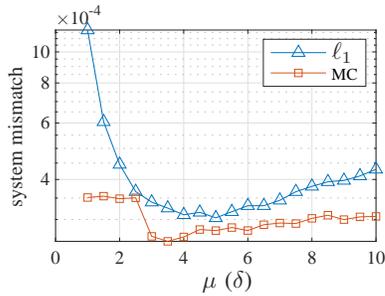}

 \caption{Performance comparisons: convex (soft shrinkage) and
 weakly-convex (firm shrinkage).}
\label{fig:l1_firm}
\end{figure}

Figure \ref{fig:l1_firm} plots the results averaged over 300 independent trials.
The approach described in Section \ref{subsec:algorithm1_with_firm}
outperforms the $\ell_1$-based approach.
Note here that it is a fair comparison in the sense that both approaches
have a single tuning parameter.
To see the full potential of the firm shrinkage operator, one may 
also tune $\lambda_1$ and $\delta$, and this will lead to further 
performance improvements at the cost of additional parameter tuning.


\section{Concluding Remarks}\label{sec:conclusion}

This paper studied
the operator-regularization approaches 
from the ``optimality'' perspective
based on the ``single-valued'' proximity operator of nonconvex
functions.
The single-valuedness here was 
a necessary condition for the operator to be continuous.
(Proposition \ref{proposition:single_valuedness}).
The fundamental result was presented, showing
that a given operator is the ``single-valued'' proximity operator of 
a weakly convex function if and only if 
it is a MoL-Grad denoiser,
i.e., the gradient of a $\kappa$-smooth convex function for $\kappa>1$
(Theorem \ref{theorem:weaklyconvex_necsuffcondition}).
Extensions of Moreau's decomposition 
were also presented, expressing a vector
in terms of the proximity operators of a weakly convex
function and the conjugate of its convexified one (Propositions
\ref{proposition:extended_Moreau_decomposition} and
\ref{proposition:moreau_decomp_mod}).
In addition to the forward-backward splitting algorithm, 
the primal-dual splitting algorithm 
employing a MoL-Grad denoiser (under the slight modification)
was shown to generate a vector sequence convergent weakly to
a minimizer of the cost function involving the implicit regularizer
(Theorems \ref{theorem:convergence} and \ref{theorem:primal_dual}).
The simulation results matched our theoretical findings
as well as demonstrating potential advantages of the proposed approach.

The proposed framework has two remarkable advantages over
the previous works on convergence analysis of the PnP method.
First, it is free from Lipschitz-constant control
during the learning of the denoiser.
Second, the implicit regularizer induced by the MoL-Grad denoiser
is not confined to smooth ones, meaning that 
our framework encompasses a large class of denoisers inducing
possibly nonsmooth regularizers.
An important implication of our findings is the following:
when the fidelity term is strongly convex and
when such a modified operator splitting algorithm
that has guarantee of convergence to an optimal point is desired,
one can restrict attention to MoL-Grad denoisers, or, equivalently,
(the proximity operator of) weakly convex functions.
Applications of the presented theory to various problems
are left as an exciting open issue.

\appendices


\newcounter{appnum}
\setcounter{appnum}{1}

\setcounter{theorem}{0}
\renewcommand{\thetheorem}{\Alph{appnum}.\arabic{theorem}}

\setcounter{proposition}{0}
\renewcommand{\theproposition}{\Alph{appnum}.\arabic{proposition}}

\setcounter{lemma}{0}
\renewcommand{\thelemma}{\Alph{appnum}.\arabic{lemma}}
\setcounter{example}{0}
\renewcommand{\theexample}{\Alph{appnum}.\arabic{example}}
\setcounter{equation}{0}
\renewcommand{\theequation}{\Alph{appnum}.\arabic{equation}}
\setcounter{claim}{0}
\renewcommand{\theclaim}{\Alph{appnum}.\arabic{claim}}
\setcounter{remark}{0}
\renewcommand{\theremark}{\Alph{appnum}.\arabic{remark}}

\section{Proof of Proposition \ref{proposition:single_valuedness}}
\label{subsec:proof_prop1}

We first show the monotonicity of
the ``set-valued'' operator
 $\mathbf{Prox}_f:\euclidspace\rightarrow 2^\euclidspace:x\mapsto 
\argmin_{y\in\euclidspace} J_x(y)$, where
$J_x:=f + (1/2)\norm{x-\cdot}^2$.
By virtue of Fermat's rule together with 
Lemma  \ref{lemma:prox_resolvent_nonconvex_case} given in Section \ref{subsec:wc_and_coco}, we have
\begin{equation*}
\hspace*{-2em}
 p\! \in\!  \mathbf{Prox}_f(x) \Leftrightarrow 0 \! \in\!  \partial J_x(p)
\Leftrightarrow  x \! \in\!  \partial J_0 (p)
\Leftrightarrow p \! \in\!  \big(\partial J_0 \big)^{-1} (x),
\end{equation*}
which  implies that $\mathbf{Prox}_f=\big(\partial J_0 \big)^{-1}$.
As the subdifferential $\partial J_0 $ is monotone \cite[Example 20.3]{combettes}, its inverse
$\big(\partial J_0 \big)^{-1}$ is also monotone \cite[Proposition 20.10]{combettes}.

For contradiction, suppose, for an arbitrarily fixed $x\in\euclidspace$,
that $\mathbf{Prox}_f(x)$ contains two distinct
vectors $T(x)$ and $T(x) +\delta$ for some $\delta\in \euclidspace\setminus\{0\}$.
The monotonicity of $\mathbf{Prox}_f$ suggests that, for every
 $\epsilon\in\real_{++}$, it holds that
\begin{align}
&~ \innerprod{T(x+\epsilon \delta) - (T(x) + \delta)}{(x+\epsilon\delta) -
 x} \geq  0 \nonumber\\
\Leftrightarrow &~  \innerprod{T(x+\epsilon \delta) - T(x) - \delta}{\epsilon\delta} \geq 0
\nonumber\\
\Leftrightarrow &~ \innerprod{T(x+\epsilon \delta) - T(x)}{\delta} \geq \norm{\delta}^2.\label{eq:frommonotonicity_of_proxf}
\end{align}
By the continuity of $T$ as well as that of the inner product,
\eqref{eq:frommonotonicity_of_proxf} reads
\begin{equation}
 0= \lim_{\epsilon\downarrow 0} \innerprod{T(x+\epsilon \delta) - T(x)}{\delta} \geq \norm{\delta}^2>0,
\end{equation}
which gives contradiction. Hence, $\mathbf{Prox}_f$ is a single-valued operator
 over $\euclidspace$, meaning that $J_x$ has a unique minimizer 
 for every $x\in\euclidspace$. \migip

\setcounter{appnum}{2}

\setcounter{lemma}{0}
\renewcommand{\thelemma}{\Alph{appnum}.\arabic{lemma}}

\section{Preservation of averaged nonexpansiveness}
\label{appendix:averagedness_preservation}

\begin{lemma}
\label{lemma:cocoercivity}
Let $T:\euclidspace\rightarrow \euclidspace$.
Then, for every $\alpha\in(0,1)$ and $\varrho \in\real_{++}$, the following
 hold.
\begin{enumerate}
 \item $T$ is $\alpha$-averaged nonexpansive if and only if 
$(c^{-1} T)\circ (c {\rm Id})$ is $\alpha$-averaged nonexpansive 
for every $c\in\real_{++}$.

 \item $\varrho  T$ is $\alpha$-averaged nonexpansive if and only if 
$T\circ (\varrho  {\rm Id})$ is $\alpha$-averaged nonexpansive.
In particular, $T$ is $\varrho $-cocoercive if and only if 
$T\circ (\varrho  {\rm Id})$ is firmly (1/2-averaged) nonexpansive.
\end{enumerate}
\end{lemma}
\begin{proof}
\noindent 1): It can be verified that
\begin{align}
& T= (1-\alpha) {\rm Id}+ \alpha N, ~\exists
 N:\euclidspace\rightarrow \euclidspace \mbox{ nonexpansive} \nonumber\\
\Leftrightarrow &~ c^{-1 }T\circ (c {\rm Id}) = 
 (1-\alpha){\rm Id} + \alpha c^{-1}N\circ (c {\rm Id}), \nonumber\\
 &\hspace*{8em} \exists
 N:\euclidspace\rightarrow \euclidspace \mbox{ nonexpansive.}
\label{eq:averagedness_equivalence1}
\end{align}
Here, $\tilde{N}:=c^{-1}N\circ (c {\rm Id})$ is nonexpansive if
$N$ is nonexpansive, because $\norm{\tilde{N}(x) - \tilde{N}(y)} =
 c^{-1} \norm{N(c x) - N(c y)} \leq c^{-1}\norm{c x- c y}
 =\norm{x-y}$.
The converse can be verified in an analogous way.
Hence, \eqref{eq:averagedness_equivalence1} implies the equivalence of
 the $\alpha$-averaged nonexpansiveness of $T$ and 
$c^{-1 }T\circ (c {\rm Id})$ for an arbitrary $c>0$.

\noindent 2): 
As the second claim can be immediately verified by
letting $\alpha:=1/2$ in the first one,
we prove the first claim in the following.
Assume that $\varrho T$ is $\alpha$-averaged nonexpansive.
Then, by Lemma \ref{lemma:cocoercivity}.1, 
$c^{-1}\varrho T \circ (c {\rm Id})$ is $\alpha$-averaged nonexpansive
for every $c\in\real_{++}$.
The specific choice of $c:=\varrho$  yields 
$\alpha$-averaged nonexpansiveness of $T \circ (\varrho {\rm Id})$.

To verify the reverse implication,
assume that $T\circ (\varrho {\rm Id})$ is $\alpha$-averaged nonexpansive.
Then, by Lemma \ref{lemma:cocoercivity}.1 again, 
$c^{-1} T \circ (c \varrho {\rm Id})$ is $\alpha$-averaged nonexpansive
for every $c\in\real_{++}$.
The specific choice of $c:=\varrho^{-1}$ yields
$\alpha$-averaged nonexpansiveness of $\varrho T$.
\end{proof}


\setcounter{appnum}{3}
\setcounter{lemma}{0}
\renewcommand{\thelemma}{\Alph{appnum}.\arabic{lemma}}



\section{Proof of Example \ref{example:vector_firm_shrinkage}}
\label{sec:proof_vector_firm}

We first prove the following lemma.
\begin{lemma}
\label{lemma:vector_shrinkage}
Let $\norm{\cdot}$ be a norm in a real Hilbert space $\euclidspace$, and
$h:\real_+\rightarrow \real$ be a nonincreasing function; i.e.,
$\real_+\ni a\leq b\in\real_+ \Rightarrow h(a)\leq h(b)$.
Assume that the s-prox operator of $h$ is well-defined.
Then, for every $\lambda\in\real_{++}$ and every $x\in\euclidspace$,
it holds that 
$^\lambda(h\circ \norm{\cdot})(x) = ~
^\lambda(h\circ \abs{\cdot})(\norm{x})
$ and
$\sprox_{\lambda h\circ \norm{\cdot}}(x) 
= 
\begin{cases}
\dfrac{\sprox_{\lambda h\circ \abs{\cdot}}(\norm{x})}{\norm{x}}x
 &
\mbox{if } x\neq 0
\\
0 & \mbox{if } x=0.
\end{cases}
$
\end{lemma}

\noindent {\em Proof of Lemma \ref{lemma:vector_shrinkage}:}
By definition, we have
\begin{equation}
 ~^{\lambda} (h\circ \norm{\cdot}) (x) =
\min_{y\in\euclidspace} \Big( h(\norm{y}) + \frac{1}{2\lambda}
\norm{y - x}^2\Big).
\label{eq:norm_moreauenvelope1}
\end{equation}
If $x=0$, $^\lambda(h\circ \norm{\cdot})(0) = ~
^\lambda(h\circ \abs{\cdot})(\norm{0})= h(0)$.
We now assume that $x\neq 0$.
Based on the orthogonal decomposition of 
$\euclidspace = \mathcal{M} \oplus \mathcal{M}^\perp$
with $\mathcal{M}:={\rm span} \{x\}$
and its orthogonal complement $\mathcal{M}^\perp$,
every $y$ can be represented by
$y= \xi x + y_{\perp}$ with
$\xi\in \real$ and $y_{\perp} \in \mathcal{M}^\perp$.
Thus, by the Pythagorean theorem,
the Moreau envelope in \eqref{eq:norm_moreauenvelope1}
can be represented as follows:
\begin{align}
 ~^{\lambda} (h\circ \norm{\cdot}) (x)
=&~
\min_{\xi\in\real,y_{\perp}\in \mathcal{M}^{\perp}}
h\Big( \sqrt{\norm{\xi x}^2 + \norm{y_{\perp}}^2}\Big)
\nonumber\\
&\hspace*{3.7em} 
+ \frac{1}{2\lambda} 
\left(
\norm{(\xi -1) x}^2
+ \norm{y_{\perp}}^2
\right)\nonumber\\
=&~ \min_{\xi\in\real}
h(\norm{\xi x})
+ \frac{1}{2\lambda} 
\norm{(\xi -1) x}^2
\nonumber\\
=&~
\min_{\xi\in\real}
h( \abs{\xi \norm{x}})
+ \frac{1}{2\lambda} 
(\xi\norm{x} - \norm{x})^2
\nonumber\\
=&~
 \min_{\zeta\in\real}
h( \abs{\zeta})
+ \frac{1}{2\lambda} 
(\zeta - \norm{x})^2
\nonumber\\
=&~
^\lambda(h\circ \abs{\cdot})(\norm{x}),
\label{eq:norm_moreauenvelope}
\end{align}
which proves the first claim of the lemma.
Here, the second equality of \eqref{eq:norm_moreauenvelope}
holds because
the minimum is achieved only when $y_{\perp}=0$
owing to the nonincreasingness of $h$.
The second claim  of the lemma can be verified in an analogous way.
\migip

\noindent {\em Proof of Example \ref{example:vector_firm_shrinkage}:}
We only prove the case of $x\neq 0$, as the case of $x=0$ is straightforward.
Letting $h:x\mapsto x$ and $\lambda:=\lambda_2$ 
in Lemma \ref{lemma:vector_shrinkage}
gives 
$^{\lambda_2} \norm{\cdot}(x) = ~^{\lambda_2} \abs{\cdot}(\norm{x})$.
This, together with $(\abs{\cdot})_{\lambda_2} =
\varphi_{\lambda_2}^{\rm MC}$,
verifies that 
$ (\norm{\cdot})_{\lambda_2} = \varphi_{\lambda_2}^{\rm MC}\circ \norm{\cdot}$.
Recall now that the MC function given in \eqref{eq:mcp} is 
a symmetric nondecreasing function.
Hence, letting $h = \varphi_{\lambda_2}^{\rm MC}$ and $\lambda:=\lambda_1$ 
in Lemma \ref{lemma:vector_shrinkage}
yields
$\sprox_{\lambda_1 \varphi_{\lambda_2}^{\rm MC} \circ \norm{\cdot}}(x) 
= (x/\norm{x})
\sprox_{\lambda_1 \varphi_{\lambda_2}^{\rm MC} \circ
\abs{\cdot}}(\norm{x})
= (x/\norm{x})
\sprox_{\lambda_1 \varphi_{\lambda_2}^{\rm MC}}(\norm{x})
= (x/\norm{x}) {\rm firm}_{\lambda_1,\lambda_2}(\norm{x})$, 
where the second equality is owing to 
the symmetry of $\varphi_{\lambda_2}^{\rm MC}$.
See Section \ref{subsec:examples} 
and also \cite{zhang,bayram15,selesnick} for the third equality.
\migip

\section{Proof of Example \ref{example:group_firm_shrinkage}}
\label{sec:proof_group_firm}

Owing to the group-wise separability of the group $\ell_1$ norm
$\norm{\cdot}_{{\rm g},1}$, 
the Moreau envelope can be expressed in the following separable form:
$^{\lambda_2}\norm{\cdot}_{{\rm g},1}(x)
= \min_{y\in\euclidspace} \big[\big(\sum_{i=1}^{G} \norm{y_i}_{\euclidspace_i}\big)
+ (1/(2\lambda_2))\norm{x - y}_{\euclidspace}^2\big]
= \sum_{i=1}^{G} 
\min_{y_i\in\euclidspace_i}
\big[
\norm{y_i}_{\euclidspace_i}
+ (1/(2\lambda_2)) \norm{x_i-y_i}_{\euclidspace_i}^2
\big]
= \sum_{i=1}^{G} {^{\lambda_2}\norm{\cdot}_{\euclidspace_i}} (x_i)$.
Thus, the Moreau enhanced model can also be expressed in a separable
form as follows:
$(\norm{\cdot}_{{\rm g},1})_{\lambda_2}(x) := 
\sum_{i=1}^{G} 
(\norm{\cdot}_{\euclidspace_i})_{\lambda_2} (x_i)$.
Hence, it follows that
$\sprox_{\lambda_1 (\norm{\cdot}_{{\rm g},1})_{\lambda_2}}(x)
= \argmin_{y\in\euclidspace} 
\sum_{i=1}^{G} 
[(\norm{\cdot}_{\euclidspace_i})_{\lambda_2} (y_i)+ 
(1/(2\lambda_1))
\norm{x_i - y_i}_{\euclidspace_i}^2]
= \big(
\sprox_{\lambda_1 (\norm{\cdot}_{\euclidspace_i})_{\lambda_2}}(x_i)
\big)_{i=1}^G
$,
where
$\sprox_{\lambda_1 (\norm{\cdot}_{\euclidspace_i})_{\lambda_2}}
= T^{\rm v\text{-}firm}_{\lambda_1,\lambda_2}
$
by Example \ref{example:vector_firm_shrinkage}.
\migip

	\bibliographystyle{IEEEtran}
	\bibliography{weaklyconvex}


 \begin{biography}[{\includegraphics[width=1in,height=1.5in,clip,keepaspectratio]{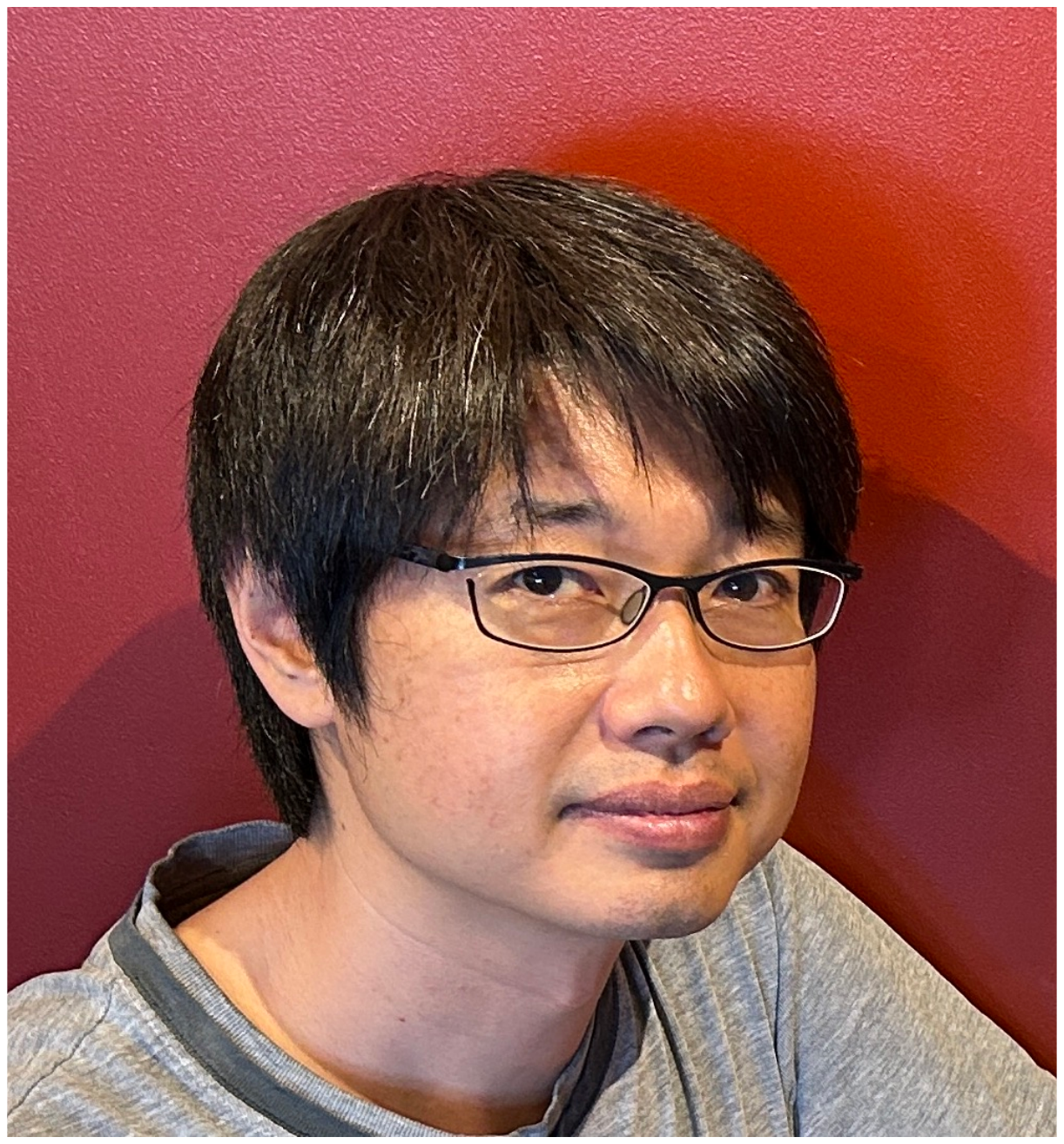}}]
{MASAHIRO YUKAWA}~(Senior Member, IEEE)~received the B.E., M.E., and
Ph.D.~degrees from the Tokyo Institute of Technology, Tokyo, Japan, in 2002, 2004,
and 2006, respectively.
He is currently a Professor with the Department of Electronics and Electrical Engineering, Keio University, Yokohama, Japan. He is the Senior Area Editor of the IEEE Transactions on Signal Processing. He served as an Associate Editor for the IEEE Transactions on Signal  Processing from 2015 to 2019, the Springer Journal of Multidimensional Systems and Signal Processing from 2012 to 2016, and the IEICE Transactions on Fundamentals of Electronics, Communications and Computer Sciences from 2009 to 2013. His research interests include mathematical adaptive signal processing, convex/sparse optimization, and machine learning.

Dr.~Yukawa received the DoCoMo Mobile Science Award in 2023, the JSPS Prize in 2022, the Young Scientists' Prize, the Commendation for Science and Technology by the Minister of Education, Culture, Sports, Science and Technology in 2014, the Excellent Paper Award from the IEICE in 2006, among many others. He is a Member of the IEICE.

 \end{biography}

 \begin{biography}[{\includegraphics[width=1.2in,height=1.8in,clip,keepaspectratio]{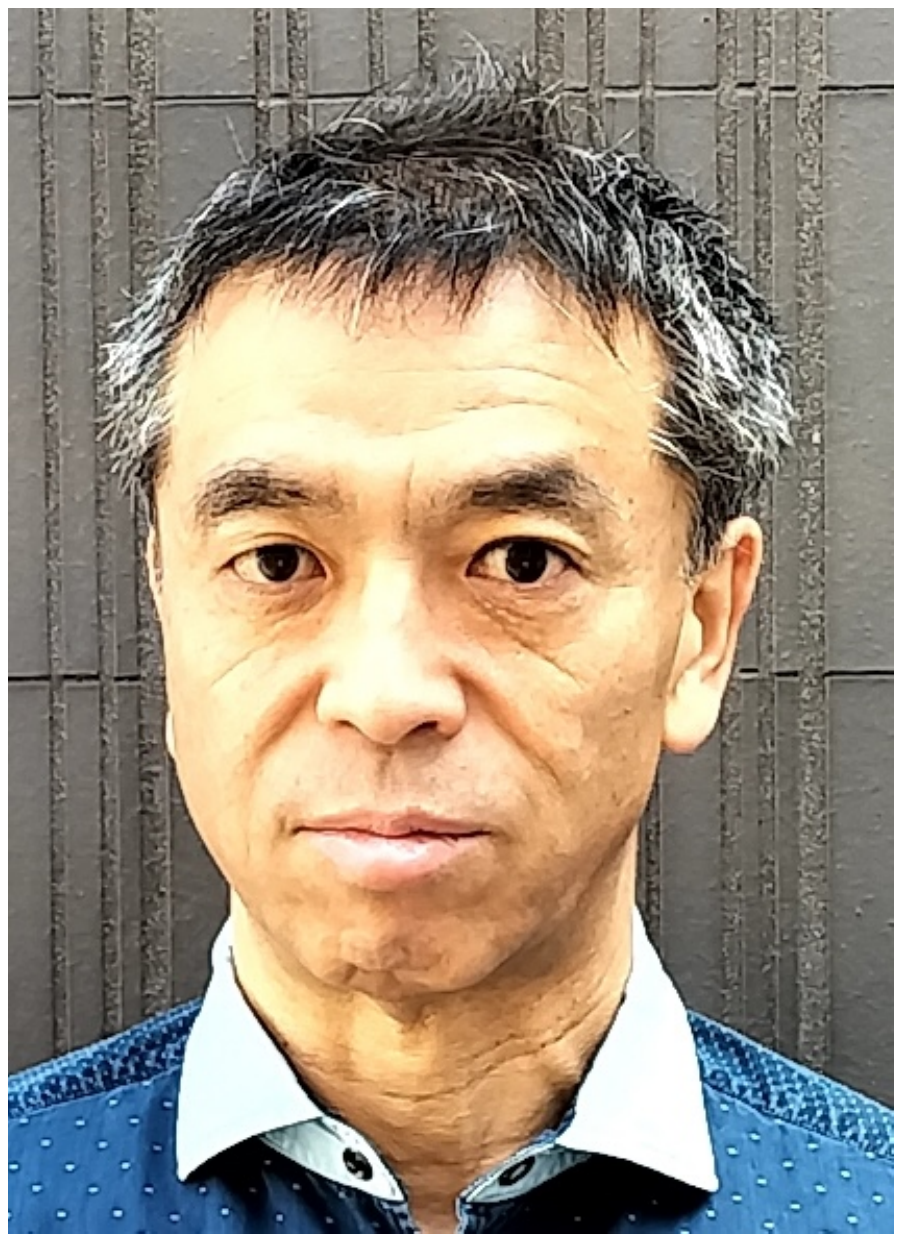}}]
{ISAO YAMADA}~(Fellow, IEEE)~received the B.E.degree in computer science 
from University of Tsukuba, in 1985, and the M.E. and Ph.D. degrees
in electrical and electronics engineering from Tokyo Institute of Technology in 1987 and in 1990
respectively. He is currently a Professor with the Department of Information and Communications
Engineering, Institute of Science Tokyo (formerly known as the Tokyo Institute of  Technology).
His research interests include mathematical signal processing, nonlinear inverse problems, and optimization theory. 
Since 2015, he has been a fellow of IEICE. He was a member of the IEEE Signal Processing Society Awards Board in 2022. 
He was a recipient of the MEXT Minister Award (Research Category), IEEE Signal Processing Magazine Best Paper Award, in 2015, 
the IEICE Excellent Paper Awards, in 1991, 1995, 2006, 2009, 2014, 2022 and 2024, 
the IEICE Achievement Award, in 2009, the ICF Research Award in 2004, 
the Docomo Mobile Science Award (Fundamental Science Division) in 2005.

 \end{biography}




%



\end{document}